\documentclass[11pt]{article}
\usepackage[english]{babel}
\usepackage[latin1]{inputenc}
\usepackage{exscale}
\usepackage[centertags]{amsmath}
\usepackage{amsfonts,amstext,amssymb,amsthm}
\usepackage{comment}
\usepackage{newlfont}
\usepackage{enumerate}
\usepackage{graphicx}
\usepackage{these}
\usepackage{color}
\usepackage{hyperref}
\usepackage{url}
\title{Generalized existence of isoperimetric regions in non-compact Riemannian manifolds and applications to the isoperimetric profile}
\author{Stefano Nardulli}
\begin{document}
      \maketitle
      \begin{center}
\noindent {\sc abstract}. For a complete noncompact Riemannian manifold with smoothly bounded geometry, we prove the existence of isoperimetric regions in a larger space obtained by adding finitely many limit manifolds at infinity. As one of many possible applications, we extend properties of the isoperimetric profile from compact manifolds to such noncompact manifolds.
\bigskip\bigskip

\bigskip

\centerline{\bf AMS subject classification: }
49Q20, 58E99, 53A10, 49Q05.
\end{center}
      \tableofcontents       
      \newpage
\section{Introduction}\label{1}
In this paper we study existence of isoperimetric regions (of minimum perimeter for prescribed volume) in a complete Riemannian manifold of infinite volume, assuming some bounded geometry conditions. The difficulty is that for a sequence of regions with perimeter approaching the infimum, some volume may disappear to infinity. We show that the region splits into a finite number of pieces which carry a positive fraction of the volume, one of them staying at finite distance and the others concentrating along divergent directions. Moreover, each of these pieces will converge to an isoperimetric region for its own volume lying in some pointed limit manifold, possibly different from the original. So isoperimetric regions exist in this generalized sense. The range of applications of these results is wide. The vague notions invoked in this introductory paragraph will be made clear and rigorous in the sequel.  
\subsection{Isoperimetric profile and existence of isoperimetric regions}
Let $(M, g)$ be a complete Riemannian manifold. We denote by $V$ the canonical Riemannian measure induced on $M$ by $g$, and by $A$ the $(n-1)$-Hausdorff measure associated to the canonical Riemannian length space metric $d$ of $M$. When it is already clear from the context, explicit mention of the metric $g$ will be suppressed in what follows. Typically in the literature, the isoperimetric profile function (or briefly, the isoperimetric profile) $I_M:[0,V(M)[\rightarrow [0,+\infty [$, is defined by $$I_M(v):= \inf\{A(\partial \Omega): \Omega\in \tau_M, V(\Omega )=v \}, v\neq 0,$$ and $I_M(0)=0$, where $\tau_M$ denotes the set of relatively compact open subsets of $M$ with smooth boundary.
However there is a more general context in which to consider this notion that will be better suited to our purposes.  Namely, we replace the set $\tau_M$ with subsets of finite perimeter, which are defined as follows. 
\begin{Def} Let $M$ be a Riemannian manifold of dimension $n$, $U\subseteq M$ an open subset, $\mathcal{X}_c(U)$ the set of smooth vector fields with compact support on $U$. Given $E\subset M$ measurable with respect to the Riemannian measure, the \textbf{perimeter of $E$ in $U$}, $ \mathcal{P}(E, U)\in [0,+\infty]$, is
      \begin{equation}
                 \mathcal{P}(E, U):=\left\{\int_{U}\chi_E div_g(X)dV_g: X\in\mathcal{X}_c(U), ||X||_{\infty}\leq 1\right\},
      \end{equation}  
where $||X||_{\infty}:=\sup\left\{|X_p|_{{g}_p}: p\in M\right\}$ and $|X_p|_{{g}_p}$ is the norm of the vector $X_p$ in the metric $g_p$ on $T_pM$. If $\mathcal{P}(E, U)<+\infty$ for every open set $U$, we call $E$ a \textbf{locally finite perimeter set}. Let us set $\mathcal{P}(E):=\mathcal{P}(E, M)$. Finally, if $\mathcal{P}(E)<+\infty$ we say that \textbf{$E$ is a set of finite perimeter}.    
\end{Def}
By standard results of the theory of sets of finite perimeter, we have that $A(\partial^* E)=\mathcal{H}^{n-1}(\partial^*E)=\mathcal{P}(E)$ where $\partial^*E$ is the reduced boundary of $E$. In particular, if $E$ has smooth boundary, then $\partial^*E=\partial E$, where $\partial E$ is the topological boundary of $E$. In the sequel we will not distinguish between the topological boundary and reduced boundary when no confusion can arise.
\begin{Def}\label{Def:IsPWeak}
Let $M$ be a Riemannian manifold of dimension $n$ (possibly with infinite volume). We denote by $\tilde{\tau}_M$ the set of  finite perimeter subsets of $M$. The function $\tilde{I}_M:[0,V(M)[\rightarrow [0,+\infty [$  defined by 
     $$\tilde{I}_M(v):= \inf\{\mathcal{P}(\Omega)=A(\partial \Omega): \Omega\in \tilde{\tau}_M, V(\Omega )=v \}$$ 
is called the \textbf{isoperimetric profile function} (or shortly the \textbf{isoperimetric profile}) of the manifold $M$. If there exist a finite perimeter set $\Omega\in\tilde{\tau}_M$ satisfying $V(\Omega)=v$, $\tilde{I}_M(V(\Omega))=A(\partial\Omega)= \mathcal{P}(\Omega)$ such an $\Omega$ will be called an \textbf{isoperimetric region}, and we say that $\tilde{I}_M(v)$ is \textbf{achieved}. 
\end{Def} 
 If $M^n$ is complete then $I_M=\tilde{I}_M$. The proof of this fact is a simple consequence of the approximation theorem of finite perimeter sets by a sequence of smooth domains, as stated in the context of Riemannian manifolds in Proposition 1.4 of \cite{MPPP}. 
\begin{Def}
We say that a sequence of finite perimeter sets $E_j$ \textbf{converge in the sense of finite perimeter sets} to another finite perimeter set $E$ if 
\begin{eqnarray*} 
\lim_{j\rightarrow+\infty} V(\Omega_j\Delta E)=0, \text{and } \lim_{j\rightarrow+\infty}\mathcal{P}(E_j)=\mathcal{P}(E).
\end{eqnarray*}
\end{Def}
For a more detailed discussion on locally finite perimeter sets and functions of bounded variation on a Riemannian manifold, one can consult \cite{MPPP}.

If $M$ is compact, classical compactness arguments of geometric measure theory combined with the direct method of the calculus of variations provide existence of isoperimetric regions in any dimension $n$. A slight variation of the latter existence argument for compact manifolds leads to the existence of isoperimetric regions in noncompact manifolds with cocompact isometry group. Roughly speaking this is because, given a minimizing sequence, isometries allow us to put the part of the volume that goes to infinity into a large geodesic ball of fixed radius, as explained in \cite{MJ}, \cite{Mor94}, or very recently in \cite{RitGalli} in the context of sub-Riemannian contact manifolds.Manuel Ritor\'e \cite{Ritsec} proved the existence of isoperimetric regions in a 2-dimensional noncompact Riemannian manifold of nonnegative sectional curvature, and very recently M Eichmair and J. Metzger \cite{EichMetz} proved the existence of isoperimetric regions for large volumes in asymptotically flat n-manifolds with nonnegative scalar curvature. Finally, if $M$ is complete, non-compact, and $V(M)<+\infty$, an easy consequence of Theorem 2.1 in \cite{RRosales} yields existence of isoperimetric regions. This is not the case for general complete infinite volume manifolds $M$. In fact it is easy to construct simple examples of $2$-dimensional Riemannian manifolds $M$ with $V(M)=+\infty$ obtained as the rotation about the $y$-axis in $\mathbb{R}^3\approx\{x,y,z\}$ of a graph in the $\{y,z\}$ plane of a smooth function, e.g.: $y\mapsto\frac{1}{y^{\alpha}}$, for some $0<\alpha\leq 1$. In this example $I_M(v)=0$ for every $v\in ]0,V(M)[$. Another example for which there are no isoperimetric regions in every volume is the familiar hyperbolic parabolid, the quadric surface $S\subset\mathbb{R}^3$ defined, for example, by the equation $z=x^2-y^2$. A deeper study of what can happen to the \textbf{subset of existence volumes of $M$} (the set of volumes $v\in ]0, V(M)[$, in which $I_M(v)$ is achieved) is given in \cite{RitCan}. In this latter example, $I_S>0$. For completeness we remind the reader that if $n\leq 7$, then the boundary $\partial\Omega$ of an isoperimetric region is smooth. More generally, the support of the boundary of an isoperimetric region is the disjoint union of a regular part $R$ and a singular part $S$. $R$ is smooth at each of its points and has constant mean curvature, while $S$ has Hausdorff-codimension at least $7$ in $\partial\Omega$. For more details on regularity theory see \cite{Morg1} or \cite{MorGMT} Sect. 8.5, Theorem 12.2.
 
In general, $I_M(v)$ is not achieved. The reason for this behavior appears in the proof of Theorem 2.1 in \cite{RRosales}, which illustrates very clearly that the lack of compactness in the variational problem is due to the fact that a part $v_1>0$ of the total volume $v=v_1+v_2$ remains at a finite distance, while another part $v_2>0$ goes to infinity. In this paper, we describe exactly what happens to the divergent part in the case of bounded geometry, which in this context means that both the Ricci curvature and the volume of geodesic balls of a fixed radius are bounded below. For any fixed $v$, we can build an example of a $2$-dimensional Riemannian manifold obtained by gluing to a Euclidean plane an infinite number (or under stronger assumptions, a finite number) of sequences of caps, diverging in a possibly infinite number of directions. Our main Theorem \ref{Main} shows that is essentially all that can occur, and that the infimum is achieved in a larger space which includes the limits at infinity.
\subsection{Main Results} 
The main result of this paper is Theorem \ref{Main}, which provides an existence result for isoperimetric regions in a noncompact Riemannian manifold satisfying the condition of bounded geometry. In the general case, solutions do not exist in the original ambient manifold, but rather in the disjoint union of a finite number of pointed limit manifolds $M_{1,\infty},\dots,M_{N,\infty}$, obtained as limit of $N$ sequences of pointed manifolds $(M, p_{ij}, g)_j$, $i\in\{1,\dots, N\}$. As a consequence of the proof of Theorem \ref{Main}, we get a decomposition lemma (Lemma \ref{algorithm}) for the thick part of a subsequence of an arbitrary minimizing sequence, which extracts exactly the structure of such sequences that is relevant to the study of the isoperimetric profile. 

Now, let us recall the basic definitions from the theory of convergence of manifolds, as exposed in \cite{Pet}. This will help us to state the main result in a precise way. 
\begin{Def} For any $m\in\mathbb{N}$, $\alpha\in ]0, 1]$, a sequence of pointed complete Riemannian manifolds is said to \textbf{converge in
the pointed $C^{m,\alpha}$ topology} (denoted $(M_i, p_i, g_i)\rightarrow (M,p,g)$) if for every $R > 0$ we can find
a domain $\Omega_R$ with $B(p,R)\subseteq\Omega_R\subseteq M$, a natural number $\nu_R\in\mathbb{N}$, and embeddings $F_{i,R}:\Omega_R\rightarrow M_i$ for large $i\geq\nu_R$ such that
$B(p_i,R)\subseteq F_{i,R} (\Omega_R)$ and $F_{i,R}^*(g_i)\rightarrow g$ on $\Omega_R$ in the $C^{m,\alpha}$ topology.
\end{Def}\noindent 
It is easy to see that this type of convergence implies pointed Gromov-Hausdorff convergence. When
all manifolds in question are closed, the maps $F_i$ are global diffeomorphisms. So for closed manifolds we can speak about unpointed convergence. What follows is the precise definition of $C^{m,\alpha}$-norm at scale $r$, that could be taken as a possible definition of bounded geometry. 
\begin{Def}[\cite{Pet}]\label{Def:BoundedGeometryPetersen}
A subset $A$ of a Riemannian n-manifold $M$  has bounded $C^{m,\alpha}$ norm on the scale of $r$, $||A||_{C^{m,\alpha},r}\leq Q$, if every point p of M lies in an open set $U$ with a chart $\psi$ from the Euclidean $r$-ball into $U$ such that
\begin{enumerate}[(i):]
      \item For all $p\in A$ there exists $U$ such that $B(p,\frac{1}{10}e^{-Q}r)\subseteq U$.
      \item $|D\psi|\leq e^Q$ on $B(0,r)$ and $|D\psi^{-1}|\leq e^Q$ on $U$.
      \item $r^{|j|+\alpha}||D^jg||_{\alpha}\leq Q$ for all multi indices j with $0\leq |j|\leq m$, where $g$ is the matrix of functions of metric coefficients in the $\psi$ coordinates regarded as a matrix on $B(0,r)$. 
\end{enumerate} 
We write that $(M,g,p)\in\mathcal{M}^{m,\alpha}(n, Q, r)$ if $||M||_{C^{m,\alpha},r}\leq Q$. 
\end{Def}

In the sequel, unless otherwise specified, we will make use of the technical assumption on $(M,g,p)\in\mathcal{M}^{m,\alpha}(n, Q, r)$ that $n\geq 2$, $r,Q>0$, $m\geq 1$, $\alpha\in ]0,1]$.  Roughly speaking, $r>0$ is a positive lower bound on the injectivity radius of $M$, i.e. $inj_M>r$, where $inj_M$ is defined below. 
\begin{Def} Let $M$ be a Riemannian manifold, set $$inj_M:=\inf_{p\in M}\{inj_{p,M}\},$$ where for every point $p\in M$, $inj_{p, M}$ is the injectivity radius at $p$ of $M$.
\end{Def} 
\begin{Def}\label{Def:BoundedGeometry}
A complete Riemannian manifold $(M, g)$, is said to have \textbf{bounded geometry} if there exists a constant $k\in\mathbb{R}$, such that $Ric_M\geq k(n-1)$ (i.e., $Ric_M\geq k(n-1)g$ in the sense of quadratic forms) and $V(B_{(M,g)}(p,1))\geq v_0$ for some positive constant $v_0$, where $B_{(M,g)}(p,r)$ is the geodesic ball (or equivalently the metric ball) of $M$ centered at $p$ and of radius $r> 0$.
\end{Def}
\paragraph{Remark:} We just mention here the fact that having Ricci tensor bounded below and positive injectivity radius implies bounded geometry in our sense.
In general, a lower bound on $Ric_M$ and on the volume of unit balls does not ensure that the pointed limit metric spaces at infinity are still manifolds, but for the proof of Theorem \ref{Main} we need the existence of limit manifolds as a hypothesis. This motivates the following definition.
\begin{Def}\label{Def:BoundedGeometryInfinity}
We say that a Riemannian manifold $(M^n, g)$ has $C^{m,\alpha}$-\textbf{bounded geometry} (at infinity) if it is of bounded geometry and if for every diverging sequence of points $(p_j)$, there exist a subsequence $(p_{{j}_{l}})$ and a pointed manifold $(M_{\infty}, g_{\infty}, p_{\infty})$ such that the sequence of pointed manifolds $(M, p_{{j}_{l}}, g)\rightarrow (M_{\infty}, g_{\infty}, p_{\infty})$, in  $C^{m,\alpha}$ topology.  
\end{Def}

\paragraph{Remark:} We observe here that Definition \ref{Def:BoundedGeometryInfinity} is weaker than Definition \ref{Def:BoundedGeometryPetersen}, and stronger than \ref{Def:BoundedGeometry}. 

In the absence of the extra condition of Definition \ref{Def:BoundedGeometryInfinity}, just assuming bounded geometry in the sense of Definition \ref{Def:BoundedGeometry}, the resulting limit space is merely a length space $(Y,y,d_Y)$ having the structure of a $C^{1,\alpha}$ manifold with a Lipschitz-Riemannian metric $g_Y$ on the regular part.  The singular part will be a set of Hausdorff dimension less than or equal to $n-1$.  
 
For more on the structure of these limit spaces, one can consult for example the works of Cheeger-Colding, \cite{ChCold1}, \cite{ChCold2}, \cite{ChCold3}. 
Regarding the smooth structure of spaces $(Y, y, d_Y)$, the reader is referred to Cheeger-Anderson \cite{AndChe92}, Anderson \cite{An92}, or \cite{Pet}, chapter 10 for a more expository discussion. In particular from \cite{AndChe92}, one can deduce that for manifolds with $Ricci$ bounded below and positive injectivity radius, the limit spaces $(Y,y,d_Y)$ are $C^{1,\alpha}$ manifolds but the convergence of the metric is just $C^{0,\alpha}$. In fact we can adjust our definition of $C^{0, \alpha}$-bounded geometry at infinity to include the set of Riemannian manifolds for which the generalized existence theorem applies. 
\paragraph{Remark:} We observe that if $|Ric_M|\leq \delta(n-1)$ and $inj_M>0$, one has by a result of M. Anderson that there exists $Q$ and $r$ such that $||(M,g)||_{C^{1,\alpha},r}\leq Q$, whose proof is contained in Theorem $76$ of \cite{Pet}. So our theorem applies in these circumstances too. 
\begin{Res}[Main: Generalized existence]
Let $M$ be a manifold of $C^{1,\alpha}$-bounded geometry. Given a positive volume $0<v < Vol(M)$, there are a finite number of limit manifolds at infinity such that their disjoint union with M contains an isoperimetric region of volume $v$ and perimeter $I_M(v)$. Moreover the number of limit manifolds is at worst linear in v.
\end{Res}
\begin{Res}[Main: Generalized existence technical statement]\label{Main}
 Let $(M,g)$ have $C^{1,\alpha}$-bounded geometry. Then for every volume $v\in]0,V(M)[$ there exist $N\in\mathbb{N}$, positive volumes $\{ v_i\}_{i\in\{1,...N\}}$, $N$ sequences of points $(p_{i,j})$,  $i\in\{1,...N\}$, $j\in\mathbb{N}$, $N$ limit manifolds $(M_{\infty,i}, g_{\infty,i}, p_{\infty, i})_{i\in\{1,...N\}}$, and domains $D_{\infty, i}\subseteq M_i$ such that  
 \begin{enumerate}[(I):]
        \item \label{Main0-} $\forall h\neq l$, $dist(p_{hj},p_{lj})\rightarrow +\infty$, as $j\rightarrow +\infty$,
        \item \label{Main0}  $(M,p_{i,j},g)\rightarrow (M_{\infty,i}, p_{\infty,i},g_{\infty,i})$ in $C^{m,\beta}$ topology for every $\beta\leq\alpha$ as $j\rightarrow +\infty$,
        \item \label{MainI}   $v=\sum_{i=1}^N v_i$,
        \item \label{MainII}  $V_{g_{\infty,i}}(D_{\infty, i})=v_i$,
        \item \label{MainIII} $A_{g_{\infty, i}}(\partial D_{\infty, i})=I_{M_{\infty,i}}(v_i)$, i.e. $I_{M_{\infty,i}}(v_i)$ is achieved by an isoperimetric region $D_{\infty,i}$,
        \item \label{MainIV}  $I_{M_{\infty,i}}(v_i)\geq I_M(v_i)$, 
        \item \label{MainIV+} if $N>1$ then $\sum_{l=1}^i I_{M_{\infty,l}}(v_i)+I_{M_{\infty, i+1}}(v-\sum_{l=1}^iv_l)\geq I_M(v)$ for every $1\leq i\leq N-1$, 
        \item \label{MainIV++} if $N>2$ then $\sum_{l=1}^i I_{M_{\infty, l}}(v_i)+I_{M_{\infty, i+1}}(v-\sum_{l=1}^iv_l)>I_M(v)$ for every $1\leq i\leq N-2$,
        \item \label{MainV}   $I_M(v)=A(\partial D_{\infty})=\sum_{i=1}^N A_{g_{\infty,i}}(\partial D_{\infty, i})=\sum_{i=1}^NI_{M_{\infty,i}}(v_i)$, where $D_{\infty}:=\bigcup_{i=1}^ND_{\infty, i}$
        \item \label{MainVI}  if $m\geq 4$ then $N\leq 2\frac{v}{v^*}+1$ and where $v^*$ depends only on $m$, $\alpha$, $n$, $Q$, $r$ is obtained in Theorem $1$ of \cite{Nar2010pv}.  
 \end{enumerate}
\end{Res}
\begin{Cor}
The conclusions $(\ref{Main0-})-(\ref{MainVI})$ are still valid under the stronger assumptions $|Ric_M|\leq k$ and $inj_M>0$, or just $Ric_M\geq k$ and $inj_M>0$.
\end{Cor}
\paragraph{Remark:} In (\ref{Main0}) if we assume $(M,g)\in\mathcal{M}^{m,\alpha}(n, Q, r)$ for some $n, Q, r$, instead of $M$ satisfying the assumption of Definition \ref{Def:BoundedGeometryInfinity}, one has to replace $\beta\leq\alpha$, with  $\beta<\alpha$. 
\paragraph{Remark:} Here $N$ depends on $v$, on the manifold $M$, and on the choice of a minimizing sequence $(D_j)_j$ in volume $v$. In fact, taking into account the example of a euclidean $2$-plane with an infinite number of bumps along $k$ diverging directions,  for every integer number $l>0$ and every volume $v$, we can construct a complete Riemannian manifold $M$ with bounded $C^{m,\alpha}$ geometry such that  every minimizing collection $\{(D_i,M_i)\}$ has exactly $l$ pieces.\\
Theorem \ref{Main} suggests the following definition.
\begin{Def} We call $D_{\infty}=\bigcup_{i} D_{\infty, i}$ a \textbf{generalized isoperimetric region of $M$}, if there exist points $p_{ij}\in M$ and pointed limit manifolds $M_{\infty, i}$, such that the $D_{\infty, i}$ are contained in the pointed limit manifolds $M_{\infty, i}$ and conditions (\ref{Main0-})-(\ref{MainIII}) and (\ref{MainV}) of the preceding theorem are satisfied.
\end{Def}
\paragraph{Remark:} We remark that $D_{\infty}$ is an isoperimetric region in volume $v$ in $\cup_i M_{\infty, i}$.
\paragraph{Remark:} If $D$ is a genuine isoperimetric region contained in $M$, then $D$ is also a generalized isoperimetric region with $N=1$ and $(M_{\infty,1}, g_{\infty,1})=(M,g)$. This does not prevent the existence of another generalized isoperimetric region of the same volume having more than one piece at infinity.
\paragraph{Remark:} The use of $C^{4,\alpha}$, or $C^{2,\alpha}$ boundedness in some subsequent arguments is due only to the technical limits of the methods employed. This assumption can be relaxed after a remark of Frank Morgan (private comunication) that avoids the use of pseudo-bubbles. Our version of that remark is contained in the second proof of the main theorem. Another way to encompass $C^{4,\alpha}$ assumptions is to notice that with a slight modification of its proof, Theorem 1 of \cite{Nar} (restated here in Lemma \ref{Regthm}) remains true even if one requires only convergence in $C^{2}$ topology of the family of metrics $g_j$. In any case the pseudo-bubble approach is much more complicated. The only delicate point that needs an explanation here is how small bounded isoperimetric regions actually can be regarded as pseudo-bubbles under the hypothesis $m=2$. What we gain using pseudo bubbles is a demonstration that the constant $v^*$ depends only on the bounds of the geometry of $M$.\\ \indent
\paragraph{Remark:} The isoperimetric regions $D_{\infty, i}$, could be disconnected and our proof works with unbounded $D_{\infty, i}$, too. Actually, the $D_{\infty, i}$ are bounded provided we have Ricci bounded below and positive injectivity radius for our manifold $M$. In case of $C^{2,\alpha}$ bounded geometry, this fact gives a simplification of the proof. Nevertheless, we prefer to give the proof without considering this simplification in order to handle our more general case. 
\paragraph{Remark:} Theorem \ref{Main} is still true if one considers minimizing cluster bubbles for a minimal partition problem instead of the isoperimetric problem. The proof is almost the same, in the sense that one has to take into account the increased complexity in notation for dealing with clusters instead that just one region. The right generalization of lemmas needed to apply the general scheme of Frank Morgan illustrated in \cite{Mor94} does not seem much more complicated than in the case of a single one region without overlapping. It could be a bit technically cumbersome to write it out in full detail, but in principle it doesn't present any major difficulty.
\paragraph{Remark:} The proof of Lemma \ref{algorithm} could be simplified by using the following theorem, whose proof is exactly the same as the Euclidean one given in Lemma 13.6 of \cite{MorGMT}.  We rewrite it in Theorem \ref{Bounded} of the Appendix for completeness. A similar use of the same Euclidean arguments is done in \cite{RRosales} Proposition 3.7.
\begin{Res}\label{Bounded} Let $M$ be a complete Riemannian manifold with bounded geometry. 
Then isoperimetric regions are bounded.
\end{Res}
Theorem \ref{Main} could be used to generalize all the results of \cite{MJ} to non-compact Riemannian manifolds with $C^{1,\alpha}$, or $C^{2,\alpha}$-bounded geometry manifolds. 
In some sense, Theorem \ref{Main} generalizes arguments which require the existence of isoperimetric regions to complete manifolds with infinite volume and $C^{1,\alpha}$, or $C^{2,\alpha}$-bounded geometry. Here we list some samples of this philosophy as corollaries of our Theorem \ref{Main}, and further applications to the existence of isoperimetric regions in complete non-compact Riemannian manifolds will appear in forthcoming papers.

\begin{CorRes}[Bavard-Pansu-Morgan-Johnson in bounded geometry]\label{BPGen}  Let $M$ have $C^{2,\alpha}$-bounded geometry, which implies that the Ricci tensor is transported to the limit manifolds. Then $I_M$ is absolutely continuous and twice differentiable almost everywhere. The left and right derivatives $I_M^-\geq I_M^+$ exist everywhere and their singular parts are non-increasing. Locally there is a constant $C_0>0$ such that $I_M(v)-C_0v^2$ is concave. Moreover, if $(n-1)K_0\leq Ricci_M$, then we have almost everywhere
\begin{equation}\label{GenBP}
          I_MI_M^{''}\leq\frac{I_M^2}{n-1}-(n-1)K_0,\text{ in the sense of distributions,}
\end{equation}
with equality in the case of the simply connected space form of constant sectional curvature $K_0\in\mathbb{R}$, (possibly $K_0\leq 0$). In this case, a generalized isoperimetric region is totally umbilic. 
\end{CorRes}
\paragraph{Remark:} Corollary \ref{BPGen} immediately extends previous results about continuity and differentiability of the isoperimetric profile that were known only for compact manifolds to the much more larger class of complete non-compact Riemannian manifolds with $C^{2,\alpha}$-bounded geometry.
\begin{CorRes}[Morgan-Johnson isoperimetric inequality in bounded geometry]\label{MJGen} Let $M$ have $C^{2,\alpha}$-bounded geometry, sectional curvature $K$ and Gauss-Bonnet-Chern integrand $G$. Suppose that
\begin{itemize}
          \item $K<K_0$, or
          \item $K\leq K_0$, and $G\leq G_0$,
\end{itemize}
where $G_0$ is the Gauss-Bonnet-Chern integrand of the model space form of constant curvature $K_0$. Then for small prescribed volume, the area of a region $R$ of volume $v$ is at least as great as $Area(\partial B_v)$, where $B_v$ is a geodesic ball of volume $v$ in the model space, with equality only if $R$ is isometric to $B_v$.   
\end{CorRes}
The proofs of Corollaries \ref{BPGen} and \ref{MJGen} run along the same lines as the corresponding proofs of theorems 3.3 and 4.4 of \cite{MJ}.
\subsection{Plan of the article}
\begin{enumerate}
           \item  Section \ref{1} constitutes the introduction of the paper. We state the main results of the paper. 
           \item In Section \ref{2}, we prove Theorem \ref{Main} as 
           stated in Section \ref{1}. 
           \item In the Appendix we prove Theorem \ref{Bounded}.
\end{enumerate}
\subsection{Acknowledgements}  
I would like to thank Pierre Pansu for actracting my attention to the subject of this paper. His comments, suggestions, and encouragements helped to shape this article. I am also indebted to Frank Morgan, Manuel Ritor\'e, Efstratios Vernadakis for their useful comments and remarks. I would like to express my gratitude to Renata Grimaldi for numerous mathematical discussions during my postdoc in Palermo were part of this manuscript was written. Finally, I thank Andrea Mondino for bringing my attention on the paper  \cite{MPPP} and Michael Deutsch for proof-reading the English of the final manuscript.    
\newpage
      \section{Proof of Theorem \ref{Main}}\label{2}
The general strategy used in calculus of variations to understand the structure of solutions of a variational problem in a noncompact ambient manifold is the Concentration-Compactness principle of P.L. Lions.  This principle suggests an investigation of regions in the manifold where volume concentrates.  For the aims of the proof, this point of view it is not strictly necessary.  But we prefer this language because it points the way to further applications of the theory developed here for more general geometric variational problems, PDE's, and the Calculus of Variations. 
\begin{Lemme}(Concentration-Compactness Lemma, \cite{Lions}  Lemma I.1)\label{C-C}
Let $M$ be a complete Riemannian manifold. Let $\mu_j$ be a sequence of Borel measures on $M$ with $\mu_j(M)\rightarrow v$. Then there is a subsequence $(\mu_j)$ such that only one of the following three conditions holds
 \begin{enumerate}[(I):]
           \item (Concentration) \label{concentration} There exists a sequence $p_j\in M$ such that for any $\frac{v}{2}>\varepsilon >0$ there is a radius $R>0$ with the property that
                                              \begin{equation}
                                                          \mu_j(B(p_j,R))>v-\varepsilon,
                                              \end{equation}
           \item (Vanishing) \label{Vanishing} For all $R>0$ there holds
                                        \begin{equation}
                                                          \lim_{j\rightarrow +\infty}Sup_{p\in M}\left\{\mu_j(B(p,R))\right\}=0,
                                        \end{equation}
           \item (Dichotomy) \label{Dichotomy} There exists a number $v_1$, $0 <v_1< v$ and a sequence of points $(p_j)$ such that for any $0<\varepsilon <\frac{v_1}{4}$ there is a number $R=R_{\varepsilon, v_1}> 0$ and two non-negative measures $\mu^1_{j,\varepsilon}$, $\mu^2_{j,\varepsilon}$ with the property that for every $R'>R$ and every strictly increasing sequence $(K_j)$ tending to $+\infty$ there exists $j_{R'}$ s.t. for all $j\geq j_{R'}$,
           \begin{eqnarray}
                     0\leq\mu_j^1+\mu_j^2\leq\mu_j ,\\
                     Supp(\mu_j^1)\subseteq B(p_j, R) for\:\: all\:\: j,\\        
                     Supp(\mu_j^2)\subseteq M-B(p_j, R'), \\
                     |\mu_j(B(p_j, R))-v_1|\leq\varepsilon,\\
                     |\mu_j^2(M)-(v-v_1)|\leq\varepsilon,\\
                     dist(Supp(\mu_j^2), Supp(\mu_j^1))\geq K_j. 
           \end{eqnarray} 
           \end{enumerate}         
\end{Lemme}
\begin{Dem}
Considering the aims of the proof, the functions of concentration $Q_j$ of Paul L\'evy are defined below. This notion serves to locate points at which volumes $v_1,...,v_N$ concentrate which are optimal in a certain sense. We define $Q_j:[0,+\infty[\rightarrow [0,v]$ by $$Q_j(R):=Sup_{p\in M}\{\mu_j(B(p,R))\}.$$ $(Q_j)$ is uniformly bounded in $BV_{loc}([0,+\infty[)$ with respect to $j$, so there exists $Q\in BV_{loc}$ such that there is a subsequence $(Q_j)$ having $j$ within $S_1\subseteq\mathbb{N}$ such that $Q_{j}\rightarrow Q(R)$ pointwise a.e. $[0,+\infty[$. Since the functions $Q_j$ are monotone increasing, so is $Q$. This ensures that the set of points of discontinuity of $Q$ is a countable set. Completing $Q$ by continuity from the left, we indeed obtain a lower semicontinuous function $Q: [0,+\infty[\rightarrow [0,v]$. It is easy to check by the theorem of existence of limit of monotone functions that there exists $v_1\in [0,v]$ such that 
\begin{equation}\label{C-C1}
           \lim_{R\rightarrow+\infty} Q(R)=v_1\in[0,v].
\end{equation}
Now, only three cases are possible: evanescence, dichotomy, and concentration. If $v_1=0$ we have evanescence, if $v_1=v$ we have concentration, and if $v_1\in ]0,v[$ we have dichotomy. Let us to explain how one can deduce (\ref{concentration})-(\ref{Dichotomy}) from (\ref{C-C1}). 
The cases $v_1=0$ and $v_1=v$ are treated exactly in the same manner as in \cite{Lions}, and we improve slightly the conclusion in the case of dichotomy.

If $v_1\in ]0,v[$ then (\ref{C-C1}) is equivalent to saying that for every $\varepsilon>0$ there is $R_{\varepsilon}>0$ such that for all $R'>R_{\varepsilon}$ 
we have 
\begin{eqnarray}
          v_1-\varepsilon<Q(R')<v_1+\varepsilon,\label{C-C1/2}\\
          v_1-\varepsilon<Q_j(R')<v_1+\varepsilon,\label{C-C1/2*}
\end{eqnarray}
for large $j$.
From (\ref{C-C1/2}) for every fixed $R>R_{\varepsilon}$ we get the existence of a sequence of points $p_{1j}$ (depending on $\varepsilon$ and $R$) with the property that
\begin{equation}\label{C-C2}
          v_1-\varepsilon<\mu_j (B(p_{1j},R))<v_1+\varepsilon,
\end{equation}
for every $j\geq j_{\varepsilon, R}$. Equation (\ref{C-C2}) is not quite what is needed for our arguments, and we must improve it to obtain exactly (\ref{Dichotomy}). This can be done by observing that if $\varepsilon$ is sufficiently small (e.g. smaller than a constant depending on $v_1$), then we can make the sequence $p_{1j}$ independent of $\varepsilon$. Following this heuristic argument, taking $\varepsilon<b_1:=\frac{1}{4}v_1$, $R_0>0$ such that $Q(R_0)>\frac{3}{4}v_1$, there exist $p_{1j}\in M$ for which (\ref{C-C2}) holds with $R$ replaced by $R_0$.  Next, take $R>0$ such that $Q(R)> v_1-\varepsilon$, so that for sufficiently large $j$ there exists a second sequence of points $p'_{1j}\in M$ for which (\ref{C-C2}) holds, and hence $$Q_j(R)+Q_j(R_0)\geq\frac{3}{4}v_1+v_1-\frac{1}{4}v_1=\frac{3}{2}v_1>\frac{5}{4}v_1>v_1+\varepsilon.$$  This implies that $B(p_{1j},R_0)\cap B(p'_{1j},R)\neq\emptyset$ for sufficiently large $j$. Thus we have $$v_1-\varepsilon<\mu_j(B(p_{1j}, R_0+2R))\leq Q_j(R_0+2R)<v_1+\varepsilon,$$ where the last inequality becomes obvious after replacing $R'$ with $R_0+2R$ in (\ref{C-C1/2*}).  This proves (\ref{Dichotomy}) with $R_{1,v_1,\varepsilon}=R_0+2R$.
\end{Dem}

This lemma will be used in our problem, taking measures $\mu_j$ having densities $\mathbf{1}_{D_j}$,
 where $\mathbf{1}_{D_j}$ is the characteristic function of $D_j$ for an almost minimizing sequence $(D_j)$ defined below. 
\begin{Def}
We say that $(D_j)_j\subseteq\tilde{\tau}_M$ (see Defn. \ref{Def:IsPWeak}) is an \textbf{almost minimizing sequence} in volume $v>0$ if \begin{enumerate}[(i):]
       \item $V(D_j)\rightarrow v$,
       \item $Area(\partial D_j)\rightarrow I_M(v)$.
\end{enumerate}
\end{Def} 
The following two Lemmas, \ref{centrale0} and \ref{centrale1}, are inspired by \cite{Mor94} Lemma 4.2 and \cite{LeoRigot} Lemma 3.1. By virtue of these we can avoid the evanescence case of Concentration-Compactness Lemma \ref{C-C}.  The difference in our treatment here is essentially in two minor changes: bounding the number of overlapping balls (which we called the multiplicity $m$ of the covering used in the proofs), and the Riemannian relative isoperimetric inequality. Both arguments use only our bounded geometry assumption, as it appears in Definition \ref{Def:BoundedGeometry}.
\begin{Lemme}\label{Lemma:Doubling}[Doubling property]\cite{EB}
Let $(M, g)$ be a complete Riemannian manifold with $Ric\geq kg$. Then for all $0<r<R$ we have 
\begin{equation}\label{Eq:Doubling}
V(B(p,R))\leq e^{\sqrt{(n-1)|k|}R}\left(\frac{R}{r}\right)^n V(B(p, r)).
\end{equation}
\end{Lemme}
\begin{Dem} The proof follows easily from the strong form of the Bishop-Gromov theorem, the fact that 
$Vol_{\lambda^2g}(B_M(x,R)))=\lambda^2Vol_g(B_M(x,\frac{R}{\lambda})),$ and the following inequalities $$\alpha_n r^{(n-1)}\leq Vol(B_{\mathbb{M}_k^n})=\alpha_n\int_0^rsinh(s)^{(n-1)}ds\leq\alpha_nr^{(n-1)}e^{r(n-1)},$$ via a conformal change of the metric. See \cite{EB}.
\end{Dem}
\begin{Cor}\label{Cor:Doubling}
Let $M^n$ be a complete Riemannian manifold with 
bounded geometry. Then for each $r>0$ there exist $c_1=c_1(n,k,r)>0$ such that $Vol(B(p,r))>c_1(n,k,r)v_0$.
\end{Cor}
\begin{Dem} If $r\geq 1$ then  $Vol(B(p,r))\geq Vol(B(p,1))\geq v_0$. If $r<1$ then (\ref{Eq:Doubling}) holds with $R=1$, hence
\begin{equation}
v_0\leq V(B(p,1))\leq e^{\sqrt{(n-1)|k|}}\left(\frac{1}{r}\right)^n V(B(p, r)).
\end{equation}
Therefore
\begin{equation}
V(B(p, r))\geq c_1(n,k,r)v_0,
\end{equation}
where $c_1(n,k,r)=Min\left\{\frac{r^n}{e^{\sqrt{(n-1)|k|}}}, 1\right\}$.
\end{Dem}
\begin{Lemme}\label{Lemma:Covering}(Covering Lemma \cite{EB}, Lemma 1.1)
Let $(M,g)$ be a complete Riemannian manifold with $Ricci\geq kg$, $k\leq 0$ and let $\rho>0$ be given. There exists a sequence of points $(x_j)\in M$ such that for any $r\geq\rho$ the following three conditions are satisfied.
\begin{enumerate}[(i):]
          \item $M\subset\cup_j B(x_j, r)$;
          \item $B(x_j, \frac{\rho}{2})\cap B(x_i, \frac{\rho}{2})=\emptyset$;
          \item if $N\in\mathbb{N}\setminus \{0\}$ is defined as $$N:=Max\{l| \exists\; j_1,...,j_l, \text{satisfying}\; \cap_{i=1}^lB(x_{j_i}, \frac{\rho}{2})\neq\emptyset\},$$ then there exists a constant $N_1=N_1(n,k, \rho, r)>0$, s.t. $N\leq N_1$.
\end{enumerate}

\end{Lemme}
\begin{Dem}
See Hebey \cite{EB} Lemma 1.1
\end{Dem}
\begin{Lemme}\label{Lemma:Relative}(Relative isoperimetric inequality \cite{MS} corollary 1.2) Let $M^n$ be a connected complete Riemannian manifold with $Ricci\geq kg$, $k\leq 0$. Then for every geodesic ball $B=B(p,r)$ and domain $D\subset M$ with smooth boundary $\partial D$ such that $V(D\cap B)\leq\frac{Vol(B)}{2}$, there exists $c(n)$ depending only on $n$ such that
\begin{equation}
V(D\cap B)^{\frac{n-1}{n}}\leq e^{c(n)(1+\sqrt{|k|}r)} rV(B(p,r))^{-\frac{1}{n}}A(\partial D\cap B).
\end{equation}
for all $p\in M$ and $r>0$. In particular, if $r=1$ and $M$ has bounded geometry, then
\begin{equation}
V(D\cap B)^{\frac{n-1}{n}}\leq e^{c(n)(1+\sqrt{|k|})} v_0^{-\frac{1}{n}}A(\partial D\cap B).
\end{equation}
\end{Lemme}
\begin{Lemme}\label{centrale0}(Non-evanescence I)
Let $M^n$ be a Riemannian manifold with bounded geometry (Defn. \ref{Def:BoundedGeometry}). Given a radius $r>0$, there are positive constants $c_3=c_3(n, k, v_0,r)>0$ and $w=w(n,k, v_0,r)>0$ such that for any set $D$ of finite perimeter $A$ and of volume $v$, there is a point $p\in M$ such that  
\begin{equation}
              V(B(p,r)\cap D)\geq m'_0=m'_0(n,k,r,v_0,v),
\end{equation} 
where $m'_0:=Min\left\{ w, c_3\frac{v^n}{A^n}\right\}.$
\end{Lemme}
\begin{Dem} Fix $r>0$. If for some point $p\in M$ one has $$V(D\cap B(x,r))\geq \frac{1}{2} V(B(p,r))\geq c(n,k,r)v_0,$$ with $c(n,k,r)=\frac{1}{2}c_1$ and $c_1$ given in Corollary \ref{Cor:Doubling}, then we can take $m'_0= w(n, k, v_0,r)=c(n, k, r)v_0$.  So assume that for all points $p$ in $M$, one has $$V(D\cap B(x,r))< \frac{1}{2} V(B(x,r)).$$ Let $\mathcal{A}$ be a maximal family of points in $M$ such that $d(x,x')\geq\frac{r}{2}$ for all $x,x'\in\mathcal{A}$ with $x\neq x'$, and $V(D\cap B(x,\frac{r}{2}))>0$ for all $x\in\mathcal{A}$. Then $V(D-\bigcup_{x\in\mathcal{A}}B(x,r)))=0$, since otherwise there would exist a point $y\in M$ such that 
\begin{equation}\label{LR1}
V\left( \left( D-\bigcup_{x\in\mathcal{A}}B(x,r)\right) \cap B(y,\frac{r}{2})\right)>0,
\end{equation}
and maximality of $\mathcal{A}$ would imply that $y\in B(x,\frac{r}{2})$ for some $x\in\mathcal{A}$, so $B(y,\frac{r}{2})\subset B(x,r)$, which contradicts (\ref{LR1}). Putting $$C=Max_{x\in\mathcal{A}}\{V(D\cap B(x,r))^{\frac{1}{n}} \},$$ we have
\begin{equation}
           \begin{array}{lll}
          V(D) & \leq & \sum_{x\in\mathcal{A}} V(D\cap B(x,r))\\
                     & \leq & C\sum_{x\in\mathcal{A}} V(D\cap B(x,r))^{\frac{n-1}{n}} .   
           \end{array}                 
\end{equation} 
A relative isoperimetric inequality for balls of radius $r$ in a Riemannian manifold with bounded geometry as stated in Lemma \ref{Lemma:Relative} (see \cite{MS} corollary 1.2) will give a constant $\gamma=\gamma(n, k, v_0,r)>0$ such that
               \begin{equation}
                          V(D\cap B(x,r))^{\frac{n-1}{n}}\leq\gamma A((\partial D)\cap B(x,r) ).  
               \end{equation}  
               Now
\begin{equation}
           \begin{array}{lll}                   
                    V(D) & \leq &  C\sum_{x\in\mathcal{A}} \gamma A((\partial D)\cap B(x, r) )\\
                      & \leq & C\gamma m A(\partial D) .
            \end{array}  
\end{equation}
where $m$ is a constant which bounds the multiplicity of the current\\ $\sum_{x\in\mathcal{A}} ((\partial D)\cap B(x,r))$ and which depends only the ratio of the volumes of balls of radii $2r$ and $\frac{1}{4}r$ of the comparison manifolds. We can take $m$ equal to $N_1(n,k,r,2r)$, where $N_1$ is the constant computed in Lemma \ref{Lemma:Covering}. Another estimate for this number is given as follows. Setting $\mathcal{A}(z):=\left\{x\in\mathcal{A} | \; z\in B(x,r) \right\}$, we observe that the balls of the family $\left\{ B(x,\frac{1}{4}r)\right\} _{x\in\mathcal{A}}$ are disjoint, and moreover
$$\cup_{x\in\mathcal{A}(z)} B(x, \frac{1}{4}r)\subseteq  B(z, 2r),$$
\begin{eqnarray*} 
card(\mathcal{A}(z))v_0 & \leq & card(\mathcal{A}(z))V(B(x,\frac{1}{4}r))
\\
& \leq & V(\cup_{x\in\mathcal{A}(z)} B(x, \frac{1}{4}r))\\
                           & \leq & V(B(z,2r))\\
                           & \leq & V(B_{k}(2r)),
\end{eqnarray*}                           
and finally
$$card(\mathcal{A}(z))\leq\frac{V(B_{k}(2r))}{v_0}=c_2(n, k, r, v_0).$$ The last inequality is a straightforward application of Bishop-Gromov's theorem. Setting $m(z):=card(\mathcal{A}(z))$, the function $m:z\rightarrow m(z)$ is exactly the multiplicity of the current  $\sum_{x\in\mathcal{A}} (\partial D)\cap B(x,r)$, and for this reason $$\sum_{x\in\mathcal{A}} A\left((\partial D)\cap B(x,r)\right)=\int_{\partial D} m(z)dV_g (z)\leq c_2(n,k, v_0, r)A(\partial D).$$ 
It follows that for some $p\in M$, 
\begin{equation}\label{imp}
V(D\cap B(p,r))^{\frac{1}{n}}\geq\frac{V(D)}{c_2(n, k, v_0, r) A(\partial D )}>0.
\end{equation} 
Taking
\begin{equation} 
m'_0:=Min\{\frac{V(D)^n}{c_2(n,k,v_0,r)^n A(\partial D )^n+1}, w(n, k, v_0,r)\},
\end{equation}
the theorem is proved.
\end{Dem}

The lemma that follows is used to avoid the evanescence case among the concentration-compactness principle alternatives for the isoperimetric problem in bounded geometry. 

\begin{Lemme}\label{centrale1}(Non-evanescence II)
Let $M$ be a Riemannian manifold with bounded geometry. Given a positive radius $r>0$, then there exist two constants $c_4=c_4(n, k, v_0,r)>0$ and $w=w(n, k, v_0,r)>0$ with the following properties. Assume that $D_j$ are currents such that $V(D_j)\rightarrow v>0$ and $I_M (v)=\lim_{j\rightarrow +\infty} Area (\partial D_j)$. Then there exists a sequence of points $p_j\in M$ such that 
\begin{equation}
              V(B(p_j,r)\cap D_j)\geq m_1=m_1(n, k, v_0, r,v).
\end{equation} 
Moreover there is an $m_2\leq m_1$ that can be choosen such that\\ $v\mapsto m_2(n,k,r,v_0, v)$ is continuous. 
\end{Lemme}

\textbf{Remark:} The importance of this lemma is that the constant $m_1$ it produces is independent of the minimizing sequence. This will guarantee that all the mass is recovered according to Morgan's scheme of the proof of the main theorem in \cite{Mor94}. The use of this lemma in the proof of the main theorem is done taking a fixed radius $r=1$. The centers $p_j$ depend on the minimizing sequence, but thankfully this is irrelevant for what follows.
 
\begin{Dem} For sufficiently large $j$, \ref{imp} implies that
            \begin{eqnarray}
                     V(D_j\cap B(p,r))^{\frac{1}{n}}\geq\frac{V(D_j)}{c_2 A(\partial D_j )} & \geq & \frac{V(D_j)}{2c_2 I_M (v)+1}\\ \nonumber
                       & \geq &\frac{v}{4c_2 A( \partial B_{M_k ^n}(\rho))+1}>0.    
            \end{eqnarray}
Here  $B_{\mathbb{M}_k ^n}(\rho_j)$ is the ball of volume $v$ in the simply connected comparison space $\mathbb{M}_k^n$ of constant sectional curvature $k$. 
We recall here that  $A( \partial B_{M_k ^n}(\rho))=A(n,k,v)$ depends only on $n, k, v$, and $v\mapsto A(n,k,v)$ for each $n,k$ fixed is continuous. Put  
\begin{equation}
m_2:=Min\{\frac{v^n}{4^nc_2(n,k ,v_0, r)^nA( \partial B_{M_k ^n}(\rho))^n+1}, w(n,k ,v_0, r)\}.
\end{equation} 
Although $v\mapsto m_1(n,k,r,v_0,v)$ is not necessarily continuous, where
\begin{equation}
m_1(n,k,r,v_0,v):=Min\{\frac{v^n}{4^nc_2^n I_M (v)^n+1}, w(n,k ,v_0, r)\},
\end{equation} 
we observe the crucial fact that $v\mapsto m_2(n,k,r,v_0,v)$ is continuous.
\end{Dem}

\subsubsection{Existence of a minimizer in a $C^{m,\alpha}$ limit manifold}
\textbf{Some known preliminary results}

Here we want to apply the theory of convergence of manifolds to the isoperimetric problem when there is a lack of compactness due to a divergence to infinity of a non-neglegible part of volume in a almost minimizing sequence.

Let us recall the basic compactness result from the theory of convergence of manifolds, as exposed in \cite{Pet}. 
\begin{Thm}\label{Fthmct}(Fundamental Theorem of Convergence Theory) [\cite{Pet} theorem $72$]. For given $Q>0$, $n\geq 2$, $m\geq 0$, $\alpha\in ]0,1]$, and $r>0$, consider the class $\mathcal{M}^{m,\alpha}(n, Q, r)$ of complete, 
pointed Riemannian $n$-manifolds $(M,p,g)$ with $||(M,g)||_{C^{m,\alpha},r}\leq Q$. 
$\mathcal{M}^{m,\alpha}(n, Q, r)$ is compact in the pointed $C^{m,\beta}$ topology for all $\beta <\alpha$.
\end{Thm}
In subsequent arguments, we will need a regularity theorem in the context of variable metrics.  
\begin{Thm}\label{Regthm}\cite{Nar}
\label{T4}
Let $M^n$ be a compact Riemannian manifold, $g_j$ a sequence of Riemannian metrics of class $C^{\infty}$ that converges to a fixed metric $g_{\infty}$ in the $C^4$ topology. Assume that $B$ is a domain of  $M$ with smooth boundary $\partial B$, and $T_j$ is a sequence of currents minimizing area under volume contraints in $(M^n , g_j )$ satisfying 
\begin{equation}V_{g_{\infty}} (B\Delta T_j)\rightarrow 0.\end{equation} 
Then in normal exponential coordinates $\partial T_j$ is the graph of a function  $u_j$ on $\partial B$.
Furthermore, for all $\alpha\in ]0,1[$, $u_j \in C^{2,\alpha}(\partial B)$ and $||u_j||_{C^{2,\alpha}(\partial B)}\rightarrow 0$ as $j\rightarrow +\infty$.
\end{Thm}\noindent
\textbf{Remark:} Loosely speaking, Theorem \ref{Regthm} says that if an integral rectifiable current $T$ is minimizing and sufficiently close in the flat norm to a smooth current, then $\partial T$ is also smooth and $\partial T$ can be represented as a normal graph over $\partial B$. In \cite{Nar}, the proof of the theorem includes a precise computation of the constants involved.\\

\noindent \textbf{Remark:} Theorems \ref{Fthmct} and \ref{Regthm} are the main reasons for the $C^4$ bounded geometry assumptions in this paper.\\

In the sequel we use often the following classical isoperimetric inequality due to Pierre B\'erard and Daniel Meyer. 
\begin{Thm}\label{BerMeyer}(\cite{BM} Appendix C). Let
$M^{n}$ be a smooth, complete Riemannian manifold, possibly with boundary, of bounded sectional curvature and positive injectivity radius.
Then given $0<\delta<1$, there exists $v_0>0$ such that any open set $U$ of volume
$0<v<v_0$ satisfies
\begin{equation}\label{BerMeyer1}
        Area(\partial U)\geq\delta c_n v^{\frac{n-1}{n}}.
\end{equation}
\end{Thm}\noindent      

\begin{Thm}\label{Berard-MeyerBoundedgeometry}(\cite{EB}, Lemma 3.2)
If $M$ is a smooth, complete Riemannian manifold with bounded geometry, then there exist a positive constant $c=c(n,k, v_0)$ and a volume $\bar{v}=\bar{v}(n, k, v_0)$ such that any open set $U$ with smooth boundary and satisfying $0\leq V(U)\leq\bar{v}$ also satisfies 
\begin{equation}
V(U)^{\frac{n-1}{n}}\leq c(n, k, v_0)Area(\partial U).
\end{equation}     
\end{Thm}
\paragraph{Remark:} The same conclusion of the preceding theorem can be obtained if one replaces $U$ by any finite perimeter set, using a customary approximation theorem in the sense of finite perimeter sets.
\paragraph{Remark:} The preceding theorem implies in particular that for a complete Riemannian manifold with Ricci curvature bounded below and strictly positive injectivity radius, we have $I_M(v)\sim c_n v^{\frac{n-1}{n}}=I_{\mathbb{R}^n}(v)$ as $v\rightarrow 0$.
\subsection{Proof of Theorem \ref{Main}.}
 For what follows it will be useful to give the definitions below. 
\begin{Def}
Let $\phi:M\rightarrow N$ a diffeomorphism between two Riemannian manifolds and $\varepsilon>0$. We say that $\phi$ is a \textbf{$(1+\varepsilon)$-isometry} if for every $x,y\in M$, $(1-\varepsilon)d_M(x,y)\leq d_N(\phi(x),\phi(y))\leq (1+\varepsilon)d_M(x,y)$.       
\end{Def}
For the reader's convenience, we have divided the proof of Theorem $1$ into a sequence of lemmas that in our opinion have their own inherent interest.
\begin{Lemme}\label{Isopcomparisoninfinity}
If $(M, p_j, g)\rightarrow (M_{\infty},p_{\infty}, g_{\infty})$ in $C^{k,\alpha}$ topology for $k\geq 1$, then 
\begin{equation}\label{Isopcomparisoninfinity1}
          I_{M_{\infty}}\geq I_M.
\end{equation}          
\end{Lemme}
\begin{Dem} We rewrite the proof appearing in \cite{Nar2010pv} Lemma 3.4 for the convenience of the reader. Fix $0<v<Vol(M)$. Let $D_{\infty}\subseteq M_{\infty}$ be an arbitrary domain of volume $v=Vol_{g_{\infty}}(D_{\infty})$. Put $r:=d_H(D_{\infty},p_{\infty})$, where $d_H$ denotes the Hausdorff distance. Consider the sequence $\varphi_j:B(p_{\infty},r+1)\rightarrow M$ of $(1+\varepsilon_j)$-isometries given by the convergence of pointed manifolds, for some sequence $\varepsilon_j\searrow 0$. Setting $D_j:=\varphi_{j}(D_{\infty})$ and $v_j:=Vol(D_j)$, it is easy to see that 
\begin{enumerate}[(i):]
          \item $v_j\rightarrow v$, 
          \item $Area_g(\partial D_j)\rightarrow Area_{g_{\infty}}(\partial D_{\infty})$.
\end{enumerate} 
Moreover, (i)-(ii) hold because $\varphi_j$ is a $(1+\varepsilon_j)$-isometry.

We now proceed with the proof of (\ref{Isopcomparisoninfinity1}) by contradiction. Suppose that there exist a volume $0<v<Vol(M)$ satisfying 
\begin{equation}\label{Isopcomparisoninfinity2}
          I_{M_{\infty}}(v)< I_M(v). 
\end{equation}
 Then there is a domain $D_{\infty}\subseteq M_{\infty}$ such that 
$$
I_{M_{\infty}}(v)\leq A_{g_{\infty}}(\partial D_{\infty})< I_M(v).
$$
As above we can find domains $D_j\subset M$ satisfying (i)-(ii). Unfortunately the volumes $v_j$ in general are not exactly equal to $v$. So we have to readjust the domains $D_j$ to get $v_j=v$, for every $j$, preserving the property $A_g(\partial D_j)\rightarrow A_{g_{\infty}}(\partial D_{\infty})$ as $j\rightarrow +\infty$, to get the desired contradiction. This can be done using the following construction that will be used in many places in the sequel. 
Examining the proofs of the deformation lemma of \cite{RitGalli} and the compensation lemma of \cite{Nar}, one can construct domains $D^{\infty}_j\subseteq B(p_{\infty},r+1)\subseteq M_{\infty}$ as small perturbations of $D_{\infty}$ such that   
\begin{equation}
    A_{g_{\infty}}(\partial D^{\infty}_j)\leq A_{g_{\infty}}(\partial D_{\infty})+c\tilde{v}_j,
\end{equation}
\begin{equation}
\tilde{v}_j\searrow 0,
\end{equation}
and 
\begin{equation}
Vol_g(\varphi_j(D^{\infty}_j))=v.
\end{equation}
The preceeding discussion shows the existence of bounded finite perimeter sets (in fact, smooth domains) $D_j:=\varphi_j(D^{\infty}_j)\subset M$ satisfying
\begin{equation}
Vol_g(D_j)=v,
\end{equation}
\begin{equation}
 |A_g(\partial D_j)-A_{g_{\infty}}(\partial D^{\infty}_j)|\rightarrow 0,
\end{equation}
again using the fact that $\varphi_j$ is a $(1+\varepsilon_j)$ isometry. Thus we have a sequence of domains $D_j$ of equal volume such that 
\begin{equation}
A_g(\partial D_j)\rightarrow A_{g_{\infty}}(\partial D_{\infty})<I_M(v),
\end{equation}
which is the desired contradiction.  The theorem follows from the fact that $v$ is arbitrary.
\end{Dem}
\begin{Lemme} Let $\tilde{M}:=M\bigcup_{i=1}^N M_{\infty,i}$ be a disjoint union of finitely many limit manifolds $(M_{\infty,i}, g_{\infty})=\lim_j(M,p_{i,j}, g)$. Then $I_{\tilde{M}}=I_{M}$.
\end{Lemme}
\begin{Dem}
It is a trivial to check that $I_{M}\geq I_{\tilde{M}}$. Observe that when $d_M(p_{ij}, p_{lj})\leq K$ for some constant $K>0$, we have $(M_{i,\infty}, g_{i,\infty})=(M_{l,j}, g_{l,\infty})$.  Thus we can restrict ourselves to the case of sequences diverging in different directions, i.e., $d_M(p_{ij}, p_{lj})\rightarrow+\infty$ for every $i\neq j$.  Then mimicking the proof of the preceding lemma, one can obtain $I_{M}\leq I_{\tilde{M}}$.
\end{Dem}

The proof of next lemma contains a construction of a decomposition of the $\varepsilon$-thick part of a subsequence of a minimizing sequence $D_j$  into a finite number of pieces.  These are obtained by cutting geodesic balls centered at concentration points, whose radius is determined by a coarea formula argument. The proof is inspired by \cite{RRosales}. 
\begin{Lemme}\label{algorithm}
Let $D_j\subseteq M$ be a minimizing sequence of volume $0<v<V(M)$. Suppose that there are $N\geq 1$ sequences of points $(p_{ij})_j$, $i\in\{1,\dots, N\}$, $N$ pointed limit manifolds $(M_{\infty,i},p_{i,\infty})$, and $N$ volumes $v_i$ such that 
\begin{enumerate}[(i):]
         \item  \label{algorithmi} $0<\sum_i^N v_i\leq v$ (possibly $\sum_i^N v_i<v$),
         \item \label{algorithmii} $(M, p_{ij})\rightarrow (M_{\infty,i},p_{i,\infty})$ in the $C^{m, \beta}$ topology for every $\beta<\alpha$,
         \item \label{algorithmiii} $d_M(p_{hj}, p_{lj})\rightarrow +\infty$, for every $h\neq l$,
         \item \label{algorithmiv} for every  $\varepsilon>0$ there exists $R_{\varepsilon}$ such that for all $R\geq R_{\varepsilon}$, there is a $j_{\varepsilon,R}$ satisfying $V(D_j\cap B(p_{ij}, R_{\varepsilon}))\in [v_i-\varepsilon, v_i+\varepsilon]$ for all $j\geq j_{\varepsilon,R}$. 
\end{enumerate}  
Then
\begin{enumerate}[(I):]
          \item \label{algorithmI} $V(D_{\infty,i})=v_i$,
          \item \label{algorithmII} $I_{M_{\infty,i}}(v_i)=Area_{g_{i,\infty}}(\partial D_{\infty,i}),\forall i\in\{1,...,N\}$,
          \item \label{algorithmIII} $I_{M_{\infty}^{(N)}}(v_1+\cdots +v_N)=\sum_{i=1}^N Area(\partial D_{\infty,i})=Area(\partial D_{\infty}^{(N)})$,
          \item \label{algorithmIV} $I_M(v)=\sum_{i=1}^N Area(\partial D_{\infty,i})+ \liminf_{j\rightarrow +\infty}Area(\partial D''_{N,j})$,
          \item \label{algorithmIV+} $V(D''_{N,j})\rightarrow v-\sum_{i=1}^N v_i$, 
          \item \label{algorithmV} 
                           \begin{eqnarray*} 
                                     I_M(v) & \geq & \sum_i^N I_{M_{i,\infty}}(v_i)+I_M(v-\sum_i^N v_i)\\
                                                & = & I_{M_{\infty}^{(N)}}(\sum_i^N v_i)+I_M(v-\sum_i^N v_i),
                           \end{eqnarray*}
\end{enumerate}
for some $D''_{N,j}\subseteq D_j$, and $D_{\infty}^{(N)}=\bigcup_{i=1}^N D_{\infty,i}$ and $M_{\infty}^{(N)}:=\bigcup_{i=1}^N M_{i,\infty}$ are disjoint unions.       
\end{Lemme}

\textbf{Remark:} (\ref{algorithmiv}) is equivalent to (\ref{algorithmiv}'): for every  $\varepsilon>0$ there exists $R_{\varepsilon}$ such that for all $R\geq R_{\varepsilon}$ we have $v_i-\varepsilon\leq\liminf_{j\rightarrow +\infty}V(D_j\cap B(p_{ij}, R))\leq\limsup_{j\rightarrow +\infty}V(D_j\cap B(p_{ij}, R))\leq v_i+\varepsilon$.\\

\begin{Dem} 
The proof will proceed by defining domains $\tilde{D}_{ij}^c\subseteq M_{\infty,i}$ (passing to a subsequence if necessary) that are images under the diffeomorphisms $F_{ij}$ (appearing in the definition of convergence of manifolds) to suitable intersections $D'_{ij}$ of the $D_j$ with balls of radii $t_{ij}$ given by the coarea formula, ultimately to obtain a finite perimeter set $D_{\infty,i}\subseteq M_{\infty,i}$ such that $\tilde{D}_{ij}^c\rightarrow D_{\infty,i}$ in the $\mathcal{F}_{loc}(M_{\infty,i})$ topology. The rate of the truncation with respect to variable $j$ must have the correct dependence on the rate of convergence of the pointed manifolds in order to get a division of a subsequence of the domains $D_j$ into disjoint geodesic balls centered at points $p_{ij}$. This will be achieved by first taking an exhaustion of $M_{\infty,i}$ by geodesic balls $B_M(p_{i,\infty}, R_{ik})$ (i.e. $M_{\infty,i}=\bigcup_{k=1}^N B_M(p_{i,\infty}, R_{ik})$) while keeping $R_{ik}$ fixed, then pulling back the domains $D_j$ on the manifold $\bigcup_{i=1}^N M_{\infty,i}$ via the diffeomorphisms $F_{ij}$. A standard compactness argument from geometric measure theory can then be applied for each of the domains inside these balls in the limit manifolds, and we obtain the desired domains $D_{\infty,i}$ by a diagonal argument. Finally, we pass to the limit and show that the latter procedure actually satisfies all the required properties (\ref{algorithmI})-(\ref{algorithmIV}).

To begin, we take a sequence of radii $(r_k)$ satisfying $r_{k+1}\geq r_k+2k$, and consider an exhaustion of $M_{\infty,i}$ by balls with centers $p_{\infty,i}$ and radii $r_k$, so that $M_{\infty,i}=\bigcup_k B(p_{\infty,i}, r_k)$. Then for every $k$, the convergence in the $C^{m,\beta}$ topology gives the existence of $\nu_{i,r_k}>0$ and diffeomorphisms $F_{ij,r_k}:B(p_{\infty,i},r_k)\rightarrow B(p_{ij},r_k)$ for all $j\geq\nu_{i,r_k}$.   The $F_{ij,r_k}$ are $(1+\varepsilon_{ij})$-isometries for the $N$ sequences $(\varepsilon_{ij})_j$, with the property that for each $i\in\{1,\dots,N\}$, we have $0\leq\varepsilon_{ij}\rightarrow 0$ as $j\rightarrow +\infty$. Put $\nu_{r_k}=Max_i\{\nu_{i,r_k}\}$.

At this stage we start the diagonal process, determining a suitable double sequence of cutting radii $t_{i,j,k}>0$ and $j\in S_k\subseteq\mathbb{N}$ for some sequence of infinite sets $\mathbb{N}\supset S_1\supseteq ...\supseteq S_{k-1}\supseteq S_k\supseteq S_{k+1}\supseteq ...$, defined inductively. Roughly speaking, we perform $N$ cuts by disjoint balls at each step $k$.
Before proceeding, we recall the argument of coarea used in this proof repeatedly. For every domain $D\subseteq M$, every point $p\in M$, and interval $J\subseteq\mathbb{R}$, there exists a $t\in J$ such that
\begin{equation}
  Area(D\cap(\partial B(p, t)))=\frac{1}{|J|}\int_J Area((\partial B(p, s))\cap D)ds\leq\frac{V(D)}{|J|}.
\end{equation}
We now proceed as follows: cutting with radii $t_{i,j,1}\in ]r_1, r_{1}+j[$ for $j\geq\nu_{r_2}$, we obtain domains $D'_{1,j}=D_j\cap B(p_j,t_{1,j})$, $D''_{1,j}=D_j-D'_{1,j}$ such that for $j$ large enough (e.g., $j\geq\tilde{\nu}_{1}=Max\{\nu^*_1,\nu_{r_2}\}$), we have 
\begin{equation}
           d(p_{hj},p_{lj})\geq\sum_{i=1}^N R_{i,v_i}+4,
\end{equation} 
for every $h\neq l$. Note that $\nu^*_{1}$ exists because of (\ref{algorithmiii}) of this lemma, and that
\begin{equation}\label{elmpv1}
\left|Area(\partial D'_{1,j})+Area(\partial D''_{1,j})-Area(\partial D_j)\right|\leq v.
\end{equation}

Consider the sequence of domains $\left(\tilde{D}_{1,j}=F_{j,r_2}^{-1}(D'_{1,j})\right)_j$ for $j\geq\tilde{\nu}_1$.  Then 
\begin{enumerate}
      \item $Area(\partial D'_{1,j})\leq Area(\partial D_j)+2\frac{v}{1}\leq I_M(v)+1+2v$,
      \item $V(D'_{1,j})\leq v$, 
\end{enumerate}
so the volume and boundary area of the sequence of domains is bounded. A standard argument of geometric measure theory allows us to extract a subsequence $D'_{1,j}$ with $j\in S_1\subseteq\mathbb{N}$ converging on $B(p_{\infty},r_2)$ to a domain $D_{\infty,1}$ in $\mathcal{F}_{B(p_{\infty},r_2)}$. Now look at the subsequence $D_j$ with $j\in S_1$ and repeat the preceding argument to obtain radii $t_{2,j}\in]r_2, r_2+2[$ and a subsequence $D'_{2,j}=D_{j}\cap B(p_j,t_{2,j})$ for $j\in S_1$ and $j\geq\nu_{r_3}$ such that 
\begin{equation}\label{elmpv2}
\left|Area(\partial D'_{i,2,j})+Area(\partial D''_{*,2,j})-Area(\partial D_{j})\right|\leq\frac{v}{2}.
\end{equation}
Again, by virtue of assumption (\ref{algorithmiii}), there exist $\nu^*_{2}$ such that for all $j\geq\nu^*_{2}$ , we have
\begin{equation}
           d(p_{hj},p_{h'j})\geq\sum_{i=1}^N R_{i,v_i}+8, h\neq h', h,h'\in\{1,...,N\}.
\end{equation}
This ensures that we can cut with pairwise disjoints geodesic balls. Analogously, the sequence $\left(\tilde{D}_{2,j}=F_{j,r_3}^{-1}(D'_{2,j})\right)_j$ for $j$ running in $S_1$ and $j\geq\tilde{\nu}_2=Max\{\nu^*_1,\nu_{r_3}\}$ has both bounded volume and bounded boundary area, so there is a convergent subsequence $\left(\tilde{D}_{2,j}\right)$ defined on some subset $S_2\subseteq S_1$ that converges on $B(p_{\infty},r_3)$ to a domain $D_{\infty,2}$ in $\mathcal{F}_{B(p_{\infty},r_3)}$. Continuing in this way, we obtain $S_1\supseteq ...\supseteq S_{l-1}\supseteq S_l$, radii
$t_{i,j,l}\in ]r_l, r_l+l[$, and domains $D'_{i,j,l}=D_{j}\cap B(p_{ij}, t_{i,j,l})$, $D''_{i,j,l}=D_{j}-D'_{i,j,l}$ satisfying 
\begin{equation}\label{elmpv3}
\left|\sum_{i=1}^NArea(\partial D'_{i,j,l})+Area(\partial D''_{*,j,l})-Area(\partial D_j)\right|\leq N\frac{v}{l},
\end{equation}
\begin{equation}
         d(p_{hj},p_{h'j})\geq\sum_{i=1}^N R_{i,v_i}+4k, h\neq h', h,h'\in\{1,...,N\}
\end{equation}
for all $1\leq l\leq k$, $j\in S_l$, $Inf(S_l)\geq\nu^*_l$ (i.e. for large $j$) and for all $k\geq 1$. Moreover, putting $\tilde{D}_{i,j,l}=F_{i,j,r_{l+1}}^{-1}(D'_{i,j,l})$, for all $1\leq l\leq k$ and $j\in S_l$ we have convergence of $(\tilde{D}_{i,j,l})_{j\in S_l}$ on $B(p_{\infty},r_{l+1})$ to a domain $D_{\infty,i,l}$ in $\mathcal{F}_{B(p_{\infty},r_{l+1})}$ for all $i\geq 1$, $1\leq l\leq k$. 
Let $j_k$ be chosen inductively so that 
\begin{eqnarray}\label{elmpv4bis}
       j_k<j_{k+1},\\
       V(\tilde{D}_{i,\sigma_{k}(j_k),k}\Delta D_{\infty,i,k})\leq\frac{1}{k},      
\end{eqnarray}
and define $\sigma(k)=\sigma_{k}(j_k)$.  Then the sequence $\tilde{D}_{i,k}^c:=F_{\sigma(k),r_{k+1}}^{-1}(D'_{i,\sigma(k),k})$ converges to $D_{\infty,i}=\bigcup_k D_{\infty,i,k}$ in the $\mathcal{F}_{loc}(M_{\infty})$ topology. Now define $t_{ik}:=t_{i,\sigma(k),k}$.
Observe here that $|t_{i,k+1}-t_{i,k}|>k$. From now on, we restrict our attention to the sequences $\bar{D}_{ik}=D_{i,\sigma_k}$, $\bar{D}'_{i,k}=D'_{i,\sigma_k}$, and $\bar{D}''_{i,k}=D''_{i,\sigma_k}$, which we will simply call $D_{i,k}$, $D'_{i,k}$, and $D''_{i,k}$ by abuse of notation.  Put also $F_{i,k}=F_{i,\sigma(k),r_{k+1}}$. Rename $k$ by $j$. From this construction, we argue that by possibly passing to a subsequence one can build a minimizing sequence $D_j$ with the following properties:
\begin{enumerate}[(a):]
       \item \label{elmpv4} $\left|\left(\sum_{i=1}Area(\partial D'_{ij})\right)+Area(\partial D''_{*j})-Area(\partial D_{j})\right|\leq 2N\frac{v}{j}$,
       \item \label{elmpv5}
$\lim_{j\rightarrow +\infty} Area_g(\partial D'_{ij})=\lim_{j\rightarrow +\infty} Area_{g_{\infty,i}}(\partial \tilde{D}^c_{ij})$,
       \item \label{elmpv6} $V(\tilde{D}^c_{ij})\rightarrow V(D_{\infty,i})=v_{\infty,i}$,  
       \item \label{elmpv7} $Area(\partial D_{\infty,i})\leq\liminf Area(\partial\tilde{D}^c_{ij})$,             
       \item \label{elmpv8} $(1+\varepsilon)v_i\geq v_{\infty,i}\geq (1-\varepsilon)v_i>0$ for every $\varepsilon$ (i.e. $v_{\infty,i}=v_i$),       
       \item \label{elmpv10} $I_{M_{\infty,i}}(v_{\infty,i})=Area(\partial D_{\infty,i})$,
       \item \label{elmpv10bis} $Area(\partial D_{\infty,i})=\liminf Area(\partial\tilde{D}^c_{j})$.
\end{enumerate}
Property (\ref{elmpv4}) follows directly by the construction of the sequences $(D'_{ij})$, and (\ref{elmpv5}) is an easy consequences of the fact that the diffeomorphisms given by the $C^{1,\beta}$ convergence are $(1+\varepsilon_j)$-isometries for some sequence $0\leq\varepsilon_j\rightarrow 0$. To prove (\ref{elmpv6}), observe that
\begin{eqnarray*}
|V(\tilde{D}_{ij}^c)-V(D_{\infty,i})| & \leq & |V(\tilde{D}_{i,j}^c)-V(D_{\infty,i}\cap B(p_{\infty,i},r_{j+1}))|\\
                                                         & +     & V(D_{\infty,i}-B(p_{\infty,i},r_{j+1}))\\ 
                                     & \leq  & V((\tilde{D}_j^c\Delta D_{\infty})\cap B(p_{\infty,i},r_{j+1}))\\
                                     & +      & V(D_{\infty,i}-B(p_{\infty,i},r_{j+1})),
\end{eqnarray*}
and so $\lim_{j\rightarrow\infty} V(\tilde{D}_{ij}^c)=V(D_{\infty,i})$ by (\ref{elmpv4bis}). On the other hand, the definition of the sets $\tilde{D}_{ij}^c$ 
gives us $\{ D_{ij}^c \}\rightarrow D$ in $\mathcal{F}_{loc}(M)$. Hence $Area(\partial D_{\infty,i})\leq \liminf_{j\rightarrow\infty}Area(\partial\tilde{D}_{ij}^c)$ by the lower semicontinuity of boundary area with respect to the flat norm in $\mathcal{F}_{loc}(M)$, which proves (\ref{elmpv7}). The first inequality in (\ref{elmpv8}) is true because every $D_{\infty,i}$ is a limit in the $L^1$ norm of a sequence of currents having volume less than $v$. The second inequality follows because for large enough $k$, the radii $r_k$ are greater than $R_{\varepsilon}$, so $V(D_{\infty,i})\geq (1-\varepsilon)v_i$. Actually $v_{\infty, i}=v_i$ by virtue of (\ref{elmpv8}) for every $\varepsilon$ arbitrarily small, and by assumption (\ref{algorithmiv}). 
To show (\ref{elmpv10}) we proceed by contradiction. Suppose that there exist domains $\tilde{E}_i\in\tau_{M_{\infty,i}}$ having $V(\tilde{E}_i)=v_{\infty,i}$ and $Area(\partial\tilde{E}_i)<Area(\partial D_{\infty,i})$. Take $N$ sequences of radii $s_{ij}\in ]t_j, t_{j+1}[$ and cut $\tilde{E}_i$ by coarea obtaining $\tilde{E}_{ij}:=\tilde{E}_i\cap B(p_{\infty,i}, s_{ij})$ in such a manner that  
\begin{equation}
       Area_{g_{\infty,i}}(\tilde{E}_{ij}\cap\partial B(p_{\infty,i}, s_{ij}))\leq\frac{v_{\infty,i}}{j}.
\end{equation}
Of course, $V_{g_{\infty}}(\tilde{E}_{ij})\rightarrow v_{\infty}$ since $s_{ij}\nearrow+\infty$ for every fixed $i\in\{1,..., N\}$.  
Now fix a point $x_{0,i}\in\partial\tilde{E}_i$ and a small neighborhood $\mathcal{U}_i$ of $x_{0,i}$. For $j$ large enough, we will have $\mathcal{U}_i\subseteq B(p_{\infty,i}, r_j)$. Push forward $\tilde{E}_{ij}$ in $M$ to obtain $E_{ij}:=F_{ij}(\tilde{E}_{ij})\subseteq B(p_{ij}, r_{j+1})$ and re-adjust the volumes by slightly modifying the $E_{ij}$ in $F_{ij}(\mathcal{U})$ contained in $B(p_{ij}, t_{j+1})$.  This produces domains $E'_{ij}\subseteq B(p_{ij}, r_{j+1})$ with the properties 
\begin{equation}
      E'_{ij}\cap D''_{*j}=\emptyset, 
\end{equation}
\begin{equation}
      \sum_{i=1}^N V_g(E'_{ij}\cup D''_{*j})=v,
\end{equation}
\begin{equation}
      Area(\partial E'_{ij})\leq Area(\partial E_{ij})+c\Delta v_{ij},
\end{equation} 
with $\Delta v_{ij}=V_g(E'_{ij})-V_g(E_{ij})$.  Note that $\Delta v_{ij}\rightarrow 0$ as $j\rightarrow +\infty$ since $V(\tilde{E}_{ij})\rightarrow v_{\infty,i}$ (i.e. $V(D'_j)\rightarrow v_{\infty,i}$) and $V(D''_{*j})\rightarrow v-\sum_{i=1}^Nv_{\infty,i}$, and that $c$ is a constant independent of $j$. Defining $D^*_j:= D''_{*j}\cup\left( \bigcup_iE'_{ij}\right)$, we have
\begin{eqnarray*}
       Area(\partial D^*_j) & \leq & \sum_{i=1}^N Area(\partial E'_{ij})+Area(D''_{*j})\\
            & \leq & \sum_{i=1}^N  (1+\varepsilon_{ij})^{n-1}Area(\partial\tilde{E}_{ij})+c\Delta v_{ij}+Area(\partial D''_{*j})\\
    & \leq &\sum_{i=1}^N ( (1+\varepsilon_{ij})^{n-1}Area(\partial\tilde{E}_{ij})+\frac{v_{\infty}}{j})+c\Delta v_{ij}+Area(\partial D''_{*j}),
\end{eqnarray*} 
and hence we get
\begin{eqnarray*}
   \liminf_{j\rightarrow +\infty} Area(\partial D^*_j) & \leq & \sum_{i=1}^N Area(\partial\tilde{E}_i)+\liminf_{j\rightarrow +\infty} Area(\partial D''_{*j})\\
   & < & \sum_{i=1}^N Area(\partial D_{\infty,i})+\liminf_{j\rightarrow +\infty} Area(\partial D''_{*j})\\
   & \leq & \liminf_{j\rightarrow +\infty}\left[\sum_{i=1}^N Area(\partial D'_{ij})\right]+\liminf_{j\rightarrow +\infty} Area(\partial D''_{*j})\\
   & \leq & I_M(v). 
\end{eqnarray*}
This means that the sequence of domains $D^*_j$ does better than the original minimizing sequence $D_j$, which is a contradiction. This proves (\ref{elmpv10}), (\ref{algorithmII}), and (\ref{algorithmIII}). The proof of (\ref{elmpv10bis}) is similar.  In fact, we only have to work with $D_{\infty,i}$ instead of $\tilde{E}_i$, since the set of regular points in $\partial D_{\infty,i}\cap M_{\infty,i}$ is open and the singular set has at least Hausdorff codimension $(n-2)$. Thanks to this fact, we can use the first variation of area at constant volume for variations with support applied to small balls centered at points in the regular set of the boundary to obtain the constant $c$.  

To finish the proof, we need one last argument that gives us (\ref{algorithmIV}). 
Putting $D''_N:= D''_{*j}$, we have 
\begin{eqnarray*}
 I_M(v) & = & \liminf_{j\rightarrow +\infty}\sum_{i=1}^N Area(\partial D'_{ij})+\liminf_{j\rightarrow +\infty} Area(\partial D''_{*j})\\
   & = & \sum_{i=1}^N I_{M_{\infty,i}}(v_{\infty,i})+\liminf_{j\rightarrow +\infty} Area(\partial D''_{*j})\\
   & = & \sum_{i=1}^N I_{M_{\infty,i}}(v_i)+\liminf_{j\rightarrow +\infty} Area(\partial D''_{*j})\\
   & \geq & \sum_{i=1}^N I_{M_{\infty}}(v_i).
\end{eqnarray*} 
which gives exactly (\ref{algorithmIV}). (\ref{algorithmIV+}) is a direct consequence of definitions and (\ref{algorithmV}) now follows easily from (\ref{algorithmI})-(\ref{algorithmIV+}).
\end{Dem}
\begin{Lemme}\label{pbestimate}
For all $n,r,Q,m,\alpha$, there exist positive constants $v_4=v_4(n,r,Q,m,\alpha)$ and $C_1=C_1(n,r,Q,m,\alpha)>0$ such that for all $M\in\mathcal{M}^{m,\alpha}(n,Q,r)$ and $0<v<v_4$ where $I_M(v)$ is achieved, 
\begin{equation}\label{pbestimate1} 
       I_M(v+h)\leq I_M(v)+C_1h v^{-\frac{1}{n}},
\end{equation}
provided that $v+h<v_4$.
\end{Lemme}
\begin{Dem} Define $v_4=Min\{1, v_0, v_1, v_2\}$. Put 
$\psi_{M,p}(\tilde{v})=Area(\beta)^{\frac{n}{n-1}}$ where $\beta$ is the pseudo-bubble of $M$ centered at $p$ and enclosing volume $\tilde{v}$. Then $\tilde{v}\mapsto\psi_{M,p}(\tilde{v})$ is a $C^1$ map, and $||\psi_{M,p}||_{C^1([0,v_4])}\leq C$ uniformly with respect to $M$ and $p$ (i.e. $C=C(n,r,Q, m,\alpha)$.  This is a nontrivial consequence of the proof of the existence of pseudo-bubbles that can be found in \cite{NarAnn}. When $v+h<v_4$, $$\psi_{M,p}(v+h)\leq\psi_{M,p}(v)+Ch.$$
\begin{equation}
      \begin{array}{lll}
    I_M(v+h) & \leq & \psi_{M,p}(v+h)^{\frac{n-1}{n}}\\
             & \leq & \psi_{M,p}(v)^{\frac{n-1}{n}}\left(1+\frac{Ch}{\psi_{M,p}(v)}\right)^{\frac{n-1}{n}}\\
             & \leq & \psi_{M,p}(v)^{\frac{n-1}{n}}\left(1+\frac{n-1}{n}C'h\right)\\
             & \leq & \psi_{M,p}(v)^{\frac{n-1}{n}}+C_1h v^{-\frac{1}{n}}\\
            & = & I_M(v)+C_1h v^{-\frac{1}{n}}.  
      \end{array}
\end{equation}
\end{Dem}

Here we give the proof of Theorem \ref{Main}, which is an intrinsic generalization of \cite{Mor94} adapted to the context of convergence of manifolds.\\

\begin{Dem} The general strategy of the proof consists in starting with an arbitrary minimizing sequence and obtaining a subsequence (again called $D_j$) and a suitable decomposition of this subsequence $D_j=\bigcup _i^N D_{ij}$ for some finite number $N$ independent of $j$. 
      To achieve this, consider the concentration functions $$Q_{1j}(R):=Sup_{p\in M}\{V(D_j\cap B_M(p,R))\}.$$ Since Lemma \ref{centrale1} prevents evanescence, the concentration-compactness argument provides a concentration volume $v_1\in[m_0(n,r=1, k,v_0, v),v]$.
      Suppose the concentration of volumes occurs at points $p_j$, and take the sequence of pointed manifolds $(M, p_j)$. By our hypothesis of $C^{m,\alpha}$-bounded geometry we obtain a limit pointed manifold $(M_{\infty,1}, p_{\infty,1})$ satisfying $$(M, p_j,g)\rightarrow (M_{\infty,1}, p_{\infty,1}, g_{\infty,1})$$ in the $C^{m, \beta}$ topology for every $\beta\leq\alpha$, up to a subsequence.
       Now we are in a position to apply Lemma \ref{algorithm} with $N=1$ to the sequences $D_j$, $p_{1j}$, and the manifold $M_{\infty,1}$, to produce (passing to a subsequence if necessary) a limit domain $D_{\infty,1}\subseteq M_{\infty,1}$, and a first sequence of discarded material $D''_{1,j}$ satisfying (\ref{algorithmI})-(\ref{algorithmV}) of Lemma \ref{algorithm}. 
      If $v_1=v$, then theorem \ref{Main} is proved with $N=1$. If $v_1<v$, we iterate the procedure considering the concentration functions $$Q_{2j}(R):=Sup_{p\in M}\{V(D''_{1,j}\cap B_M(p,R))\}$$ of the sequence of discarded domains $D''_{1,j}$, and again we extract a subsequence of points $p_{2j}\in M$ (always up to a subsequence).  By the concentration-compactness lemma, they satisfy $dist(p_{1j}, p_{2j})\rightarrow +\infty$.  We now take a limit manifold $(M, p_{2,j}, g)\rightarrow (M_{\infty,2}, p_{\infty,2}, g_{\infty,2})$ using again our assumption on $C^{m,\alpha}$-bounded geometry 
      and apply Lemma \ref{algorithm} with $N=2$, and we iterate this procedure.      
    Suppose that we have carried out the construction until step $i$, so that there exist $D_{\infty,1}\subseteq M_{\infty,1}, D_{\infty,2}\subseteq M_{\infty,2}, \dots, D_{\infty,i}\subseteq M_{\infty,i}$, such that if $v_i=V_{g_{\infty,i}}(D_{\infty,i})$, $D_{\infty}^{(i)}=\bigcup_{l=1}^i D_{\infty,i}$ and $M_{\infty}^{(i)}=\bigcup_{l=1}^i M_{l,\infty}$ disjoint union,
      \begin{eqnarray}
                I_{M_{\infty,l}}(v_l)=A_{g_{l,\infty}}(\partial D_{\infty,l}),\forall l\in\{1,...,i\},\label{Main21}\\
                I_{M_{\infty}^{(i)}}(v_1+\cdots +v_i)=\sum_{l=1}^i A(\partial D_{\infty,l})=A(\partial D_{\infty}^{(i)}),\label{Main22}\\
                I_M(v)=\sum_{l=1}^i A(\partial D_{\infty,l})+ \liminf_{j\rightarrow +\infty}A(\partial D''_{i,j}),\label{Main23}\\
                I_{M_{\infty,N+1}}(v''_{N+1})\geq\liminf_{j\rightarrow +\infty}A(\partial D''_{N,j}).\label{Main24}
      \end{eqnarray}     
    If $\sum_{l=1}^i v_l=v$, the algorithm finishes and the theorem is proved with $N=i$. If $\sum_{l=1}^i v_l<v$, we continue analogously at step $i+1$, applying first the concentration-compactness lemma to the discarded materials $D''_{i,j}$ to get points $p_{i+1,j}$ and then Lemma \ref{algorithm} with $N=i+1$, and $i+1$ sequences of points $\left\{(p_{1,j}),\dots, (p_{i+1,j})\right\}$. The algorithm stops in $N$ steps if and only if all the mass $v$ is recovered (i.e. if $\sum_{i=1}^{N-1}v_i<v$ and $\sum_{i=1}^Nv_i=v$) which is equivalent to the first occurrence of the concentration case at step $N$. We must mention here that assumption (\ref{algorithmIV}) of Lemma \ref{algorithm} is satisfied because it corresponds exactly to the improved version of the concentration-compactness theorem. In fact, $V(D''_{N,j})\rightarrow v-\sum_{i=1}^{N-1}v_i$ for all $N$ as $j\rightarrow +\infty$.            
             Suppose by contradiction that (\ref{Main24}) is false.  Then there exists a domain $E\subseteq M_{\infty,N+1}$ such that $V(E)=v''_N$ and $A(\partial E)<\liminf_{j\rightarrow +\infty}A(\partial D''_{N,j})$.  We can take the improved decomposition $E\bigcup_{i=1}^N D_{\infty,i}$ to get  $V(E\bigcup_{i=1}^N D_{\infty,i})=v$ and $A(\partial (E\bigcup_{i=1}^N D_{\infty,i}))<I_M(v)$.  Now projecting $E\bigcup_{i=1}^N D_{\infty,i}$ onto the manifold $M$, we obtain an improved almost minimizing sequence of volume $v$, implying $I_M(v)<I_M(v)$, which is the desired contradiction. 
             
          The whole construction produces at most countably many domains $D_{\infty,i}$ which do not a priori recover all the volume $v$. So to finish the proof, we must show that all the volume $v$ is in fact recovered, and that the number of domains $D_{\infty,i}$ is finite.     
                      
          Proceeding by contradiction, suppose that $N\rightarrow +\infty$. First we will show that  
          \begin{equation}
          \sum_{i=1}^{+\infty}v_i=v.
          \end{equation}
          It is obvious that $v_N\searrow 0$, because the series $\sum_{i=1}^{+\infty}v_i$ is convergent. Put $$\sum_{i=1}^{+\infty}v_i=\bar{v}\leq v,$$ and 
          \begin{eqnarray}
           \tilde{m}_0=\tilde{m}_0(n,k,v_0,v-\bar{v}) & = & m_1(n,k,1, v_0, v-\bar{v})\\ \nonumber
           & = & Min\{\frac{(v-\bar{v})^n}{c_4I_M(v-\bar{v})^n+1}, w\},
           \end{eqnarray}
           where $c_4=c_4(n,k,v_0)$, $w=w(n,k,v_0)$.
           Then for large $N$ we can take $v_N$ such that 
          \begin{equation}\label{Eq:Contadiction0}
                v_N<\tilde{m}_0,
          \end{equation}
          because $v_N\rightarrow 0$. The crucial fact here is to realize that $\tilde{m}_0$ is constant w.r.t. variable $N$.
          Now suppose by contradiction that $\bar{v}<v$.  We wish to use Lemma (\ref{centrale0}) with $r=1$ applied to some $D''_{i,j}$ (the discarded material) 
          for sufficiently large $i,j$, to get a contradiction. We recall here that 
          \begin{equation}\label{recovering}
                   V(D''_{N,j})\rightarrow v-\bar{v}, \text{as}\; N\rightarrow+\infty,
          \end{equation}
          by construction. Since $Q_{N,j}(R)$ is defined as the supremum of the volume of the intersection of a domain with geodesic balls of radius $R$ over the family of all geodesic balls, one easily checks that
           \begin{equation} \label{Eq:intermediate0}
                   Q_{N,j}(R)\geq Q_{N,j}(1),\;\;\forall R\geq 1,
          \end{equation} 
          because $Q_{N,j}(.)$ is monotone nondecreasing. Next
          \begin{equation}\label{Eq:intermediate1}
          Q_{N,j}(1)\geq (V(B(p,1)\cap D''_{N-1,j}),
          \end{equation}
          for some $p\in M$ given by Lemma \ref{centrale0} applied to $D''_{N-1,j}$ taking $r=1$. Again by Lemma \ref{centrale0} we have
           \begin{eqnarray}\label{Eq:intermediate2}
           V(B(p,1)\cap D''_{N-1,j})\geq Min\{\frac{V(D''_{N-1,j})^n}{c_3(n,k,v_0)A(\partial D''_{N-1,j})^n+1}, w(n,k,v_0)\}.
           \end{eqnarray} 
           The same argument leading to (\ref{Main24}) proves that
           \begin{equation}\label{Main23bis}
                     \liminf_{j\rightarrow +\infty}A(\partial D''_{N-1,j})\leq I_M(v-\bar{v}).
           \end{equation}
           \paragraph{Remark:}
           One can show that 
            \begin{equation}
                     \liminf_{j\rightarrow +\infty}A(\partial D''_{N-1,j})=I_M(v-\bar{v}).
           \end{equation}  
           Combining  (\ref{Eq:intermediate0}), (\ref{Eq:intermediate1}), (\ref{Eq:intermediate2}), (\ref{Main23bis}), and (\ref{Main24}), and   
           taking the limit $j\rightarrow +\infty$, we have
          \begin{equation}\label{Eq:recoveringmasscontradiction} 
                  Q_{N}(R)=\lim_{j\rightarrow +\infty}Q_{N,j}(R)\geq\tilde{m}_0>0,\; \forall R\geq 1.
          \end{equation}  
        On the other hand, by the definition of $v_N$ we have
          \begin{equation}
          v_N=\lim_{R\rightarrow +\infty}Q_{N}(R).
          \end{equation}       
           Finally, taking limits over $R$ in \ref{Eq:recoveringmasscontradiction} we get
          \begin{equation}
                 v_N=\lim_{R\rightarrow +\infty}Q_{N}(R)\geq \tilde{m}_0,
          \end{equation}       
which contradicts (\ref{Eq:Contadiction0}).  So all the volume is recovered and $\sum_{i=1}^{+\infty}v_i=v$ holds.

To prove that there are only finitely many terms in the series $\sum_{i=1}^{+\infty}v_i$, there are two quite different possible arguments. The first uses the theory of pseudo-bubbles. The second follows from a private communication with Frank Morgan. Roughly speaking, both rely on the fact that in bounded geometry, the profile is asymptotically Euclidean for small volumes and that in Euclidean space, two balls cannot do better than one, which is implicit in the strict concavity of $I_{\mathbb{R}^n}$.

\textbf{First proof.} First we assume that $m\geq 4$ and prove that $N<+\infty$ by contradiction. If $N$ were unbounded, then there would exist sequences $v_{N+1}$, $v''_{N+1}$ such that 
       \begin{enumerate}[(i):]
             \item $0<v_{N+1}+v''_{N+1}=v''_N$,
             \item $I_{M_{\infty, N+1}}(v_{N+1})$ is achieved, 
             \item $I_{M_{\infty, N+1}}(v''_N)\geq I_{M_{\infty, N+1}}(v_{N+1})+I_{M}(v''_{N+1})$, 
             \item\label{item:smallvolumesplittingratio} $\frac{v''_{N+1}}{v_{N+1}}\rightarrow 0$, as $N\rightarrow +\infty$,
       \end{enumerate}
Properties (i)-(ii) are easy consequences of the proof of Lemma \ref{algorithm}. 
Property (iii) is deduced immediately from (\ref{Main24}).  Property (iv) follows from (iii) and the asymptotic form of Berard-Meyer for manifolds with bounded geometry (Lemma \ref{Berard-MeyerBoundedgeometry}), recalling that $$\lim_{j\rightarrow +\infty}V(D''_{N,j})=v''_N \text{ for every}\; N.$$ More precisely, we have\\ 
$$I_{M_{\infty, N+1}}(v_{N+1}+v''_{N+1})\sim c(n) (v_{N+1}+v''_{N+1})^{\frac{n-1}{n}},$$
$$I_{M_{\infty, N+1}}(v_{N+1})\sim c(n) v_{N+1}^{\frac{n-1}{n}},$$ 
$$I_M(v''_{N+1})\sim c(n) \left(v''_{N+1}\right)^{\frac{n-1}{n}}.$$ 
These asymptotic relations plus (iii) imply
$$(v_{N+1}+v''_{N+1})^{\frac{n-1}{n}}\sim v_{N+1}^{\frac{n-1}{n}}+(v''_{N+1})^{\frac{n-1}{n}}.$$ The latter relation is true if and only if $\frac{v''_{N+1}}{v_{N+1}}\rightarrow 0$ or $\frac{v''_{N+1}}{v_{N+1}}\rightarrow +\infty$, but the former case is excluded, because our construction gives 
\begin{equation}\label{eucl1}
         v_{N+1}\geq c(v_{N+1}+v''_{N+1}),
\end{equation}
          which is easily checked using Lemma \ref{centrale1} and the fact that Lemma \ref{BerMeyer} holds for small volumes.  Thus we can assume that $\frac{v''_{N+1}}{v_{N+1}}\rightarrow 0$.
We argue from (iii) that
\begin{equation}
\dfrac{I_{M_{\infty, N+1}}(v_{N+1}+v''_{N+1})-I_{M_{\infty, N+1}}(v_{N+1})}{v''_{N+1}}\geq\dfrac{I_M(v''_{N+1})}{v''_{N+1}}. 
\end{equation} 
Hence\footnote{Note that this is the only point at which we use the theory of pseudo-bubbles.} by (\ref{pbestimate1}),
\begin{eqnarray} 
\dfrac{I_{M_{\infty, N+1}}(v_{N+1}+v''_{N+1})-I_{M_{\infty, N+1}}(v_{N+1})}{v''_{N+1}} &\leq & C_1(v_{N+1})^{-\frac{1}{n}},\\
            \dfrac{I_M(v''_{N+1})}{v''_{N+1}}  & \leq & C_1v_{N+1}^{-\frac{1}{n}}.
\end{eqnarray}
Finally, since $\frac{I_M(v''_{N+1})}{v''_{N+1}}\sim c(n)\frac{1}{(v''_{N+1})^{\frac{1}{n}}}$, asymptotically we have
$$c(n)\frac{1}{(v''_{N+1})^{\frac{1}{n}}}\leq C_1\frac{1}{(v_{N+1})^{\frac{1}{n}}},$$ hence 
\begin{equation}\label{Main7}
          \frac{c(n)}{C_1}\leq\left(\frac{v''_{N+1}}{v_{N+1}}\right)^{\frac{1}{n}}, 
\end{equation}
which contradicts (iv). 
Thus the procedure indeed finishes in a finite number of steps, which proves (\ref{MainV}). 

\textbf{Second proof.} Here we explain our version of Morgan's suggestion. 
The argument proceeds by contradiction. For large $N$, $v''_N$ is small, and we can improve $D_{\infty}^{(N)}$ by perturbing one of its pieces in a controlled way (for example, $D_{\infty,1}$ inside $M_{\infty,1}$), obtaining a domain $\tilde{D}_{\infty,1}\subseteq M_{\infty,1}$ satisfying 
    \begin{eqnarray}
              V(\tilde{D}_{\infty,1})=v_1+v''_N,\label{Main31}\\
             Area(\partial \tilde{D}_{\infty,1})\leq Area(\partial D_{\infty,1})+c(D_{\infty,1})v''_N.\label{Main32}         
    \end{eqnarray} 
    Here $c(D_{\infty,1})$ is the integral of the mean curvature of $D_{\infty,1}$ at its regular points (see \cite{RitGalli}). Observe that the estimate (\ref{Main32}) is possible because $V(\tilde{D}_{\infty,1}\Delta D_{\infty,1})$ is small.
  Denote $\tilde{D}_{\infty}^{(N)}:=\tilde{D}_{\infty,1}\bigcup_{i=2}^N D_{\infty,i}$ for the disjoint union and observe that 
    \begin{equation}
             V(\tilde{D}_{\infty}^{(N)})=v.
    \end{equation}                 
   Combining (\ref{Main23}), (\ref{Main31}), (\ref{Main32}) yields
   \begin{eqnarray}\label{Main4}
              \sum_{i=2}^N Area(\partial D_{\infty,i}) +cv''_{N} & \geq & Area(\partial \tilde{D}_{\infty}^{(N)})\\
              & \geq & I_{\cup_iM_{i,\infty}}(v)\\
              & \geq & I_M(v)\\
              & \geq & \sum_{i=2}^N Area(\partial D_{\infty,i}) + I_M(v''_N),
   \end{eqnarray}
   where $c=c(D_{\infty,1})$.
   As an easy consequence of the preceding inequalities we have that
   \begin{equation}\label{Main5}
             c(D_{\infty,1})v''_N \geq I_M(v''_N),
   \end{equation}
   and finally 
   \begin{equation}\label{Main6}
             c(D_{\infty,1}) \geq\frac{I_M(v''_N)}{v''_N}.
   \end{equation}
   Letting $N\rightarrow +\infty$ in (\ref{Main6}), $v''_N$ becomes arbitrarily small, in particular smaller than $\bar{v}$ of Lemma \ref{Berard-MeyerBoundedgeometry}.  For small volumes, we can now apply the generalized version of Berard-Meyer in our setting (i.e., Lemma \ref{Berard-MeyerBoundedgeometry}) giving 
   \begin{equation}\label{eq:finalestimates}
   c(D_{\infty,1})\geq\frac{I_M(v''_N)}{v''_N}\geq c(n, k, v_0) \frac{(v''_N)^{\frac{n-1}{n}}}{v''_N}\geq\frac{c(n, k, v_0)}{(v''_N)^{\frac{1}{n}}}.
   \end{equation} 
   Finally, letting $v_N\rightarrow 0$ in (\ref{eq:finalestimates}), we get the desired contradiction, and thus $N$ has to be finite. 
   
   To obtain an estimate on $N$, one way to proceed is to show that there is no dichotomy for volumes less than some fixed $v^*=v^*(n, k, v_0)>0$. Assuming the existence of such a $v^*$, we observe that the algorithm produces $v_1\geq v_2\geq\cdots\geq v_N$.  Furthermore, $v_N\leq v^*$ because it is the first time that dichotomy cannot appear, which yields
\begin{equation}
v^*N\leq v_N N\leq\sum_{i=1}^N v_i=v.
\end{equation}
Consequently 
\begin{equation}\label{eq:Mainestimates}
N\leq\left[\frac{v}{v^*}\right]+1,
\end{equation}
where $v^*=v^*(n, k, v_0)$ can be taken equal to $\bar{v}$ of Lemma \ref{Berard-MeyerBoundedgeometry}.
On the other hand, one can construct examples such that for every $v$, there are exactly $N=\left[\frac{v}{v^*}\right]+1$ pieces. So in this sense, the estimate (\ref{eq:Mainestimates}) is sharp and it concludes the proof of the theorem.
\end{Dem}

Note that the argument given here furnishes only the existence of such a $v^*$, but not the optimal value.
\paragraph{Question:} What is the largest volume, i.e., what is the sharp value for $v^*$?
\paragraph{Remark:} A rigorous proof that there is no dichotomy for small volumes is contained in \cite{Nar2010pv}, via the theory of pseudo-bubbles.  This argument requires the stronger condition of $C^{4,\alpha}$-bounded geometry. What do we gain with this argument? We gain that $v^*=v^*(n, k, v_0)$ depends a priori on the bounds $Q$ on the geometry.
\paragraph{Remark:} We must observe that the second proof, even if more elegant and avoiding the use of pseudo-bubbles, introduces constants that depend on some limit domain. So for a better estimate of $N$, we have to use more involved arguments, as for example the first proof.\\ 
\paragraph{Remark:} The geometric measure theoretic compactness arguments are to be understood intrinsically, using the theory of 
sets of finite perimeter in manifolds as done in \cite{RRosales}. For the general theory one can consult  \cite{MPPP}
\section{Appendix: Proof of Theorem \ref{Bounded}}

\begin{Dem}
Fix a point $p\in M$. Let $D$ such that $I_M(v)=Area(\partial D)$ and $V(D)=v$, $V(r):=V(D\cap(M\setminus B(p,r))$, $A(r):=Area((\partial D)\cap(M\setminus B(p,r)))$. The proof proceeds exactly as in the Euclidean case but the constants involved depend on (old version the injectivity radius and) lower bounds on the Ricci curvature and volume of unit balls. As is easily seen using standard Riemannian comparison geometry techniques and Theorem \ref{Berard-MeyerBoundedgeometry},
          \begin{equation}\label{Bounded0}
                    V(r)\searrow 0
          \end{equation} as $r\rightarrow +\infty$, so the function $V(r)$ is decreasing. By the coarea formula we have $V'(r)=-Area(D\cap\partial B(p,r))$ for almost all $r$. Fix a $\delta>0$ in Lemma \ref{BerMeyer} and take $r$ large enough to have $V(r)<v_0$ where $v_0$ is obtained in Lemma \ref{BerMeyer}. From (\ref{BerMeyer1}) applied to $D\cap(M\setminus B(p,r))$ we obtain
     \begin{eqnarray}\label{Bounded1}
                |V'(r)|+A(r) & \geq & Area(D\cap\partial B(p,r))+A(r)\\
                                  & = & Area(\partial (D\cap(M\setminus B(p,r)))\\
                                  & \geq & \delta c(n)V(r)^{\frac{n-1}{n}},
     \end{eqnarray} 
     for almost all $r$. By using the fact that $D$ is an isoperimetric region, the first variation formula for area with respect to volume yields
    \begin{equation}\label{Bounded2}
           |V'(r)|+CV(r)\geq A(r).   
    \end{equation}   
    For almost all $r$ large enough we have
    \begin{equation}\label{Bounded3}
               CV(r)\leq\frac{c(n)\delta}{2}V(r)^{\frac{n-1}{n}},
    \end{equation}
    in virtue of (\ref{Bounded0}) and an application of Lemma \ref{Berard-MeyerBoundedgeometry} to $$U=D\cap(M\setminus B(p,r)),$$ since the volume of $U$ is small.
    Combining (\ref{Bounded0})-(\ref{Bounded3}) we get for large $r$
    \begin{equation}
            |V'(r)|\geq \frac{c(n)\delta}{2}V(r)^{\frac{n-1}{n}}. 
    \end{equation}
    Assuming that $D$ is unbounded implies that 
    \begin{equation}\label{Bounded3+}
             V(r)>0,
    \end{equation} for all $r$. Now, taking into account (\ref{Bounded3}) we get
    \begin{equation}\label{Bounded4}
              (n+1)\left(V(r)^{\frac{1}{n+1}}\right)'=-|V'(r)|V(r)^{-\frac{n-1}{n}}\leq -\frac{c}{4}, 
    \end{equation} 
    which contradicts (\ref{Bounded3+}).
    \end{Dem}    
      \newpage
      \markboth{References}{References}
      \bibliographystyle{alpha}
      \bibliography{these}
      \addcontentsline{toc}{section}{\numberline{}References}
      \emph{Stefano Nardulli\\ Departamento de M\'etodos Matem\'aticos\\ Instituto de Matem\'atica\\ UFRJ-Universidade Federal do Rio de Janeiro, Brasil\\ email: nardulli@im.ufrj.br} 
\end{document}